\documentclass[a4paper,12pt]{amsart}
\usepackage{a4wide,fullpage,setspace}
\usepackage{amsmath,mathtools,amssymb,abraces,enumitem}
\usepackage{amscd,graphicx,tikz}
\usepackage[all,cmtip]{xy}
\usepackage{stmaryrd,accents,microtype}
\usepackage[T1]{fontenc}
\usepackage{lmodern}
\usepackage{hyperref}
\definecolor{dark-red}{rgb}{0.5,0.15,0.15}
\hypersetup{
	colorlinks   = true, 
	urlcolor     = dark-red, 
	linkcolor    = black, 
	citecolor   = dark-red 
}
\usepackage[numbers,sort&compress]{natbib}

\title{Directed degeneracy maps for precubical sets}

\author[P. Gaucher]{Philippe Gaucher}

\address{Universit\'e Paris Cit\'e, CNRS, IRIF, F-75013, Paris, France}

\urladdr{http://www.irif.fr/{\~{}}gaucher}



\date{2023/03/09}

\swapnumbers

\renewcommand{\thepart}{\Roman{part}}
\makeatletter
\def\@part[#1]#2{%
	\ifnum \c@secnumdepth >\m@ne
	\refstepcounter{part}%
	\addcontentsline{toc}{part}{\thepart\hspace{1em}#1}%
	\else
	\addcontentsline{toc}{part}{#1}%
	\fi
	\markboth{}{}%
	{\centering
		\interlinepenalty \@M
		\normalfont
		\ifnum \c@secnumdepth >\m@ne
		\Large\bfseries \thepart\hspace{1em}
		\fi
		\Large\bfseries #2\par}
	\nobreak
	\vskip 3ex
	\@afterheading}
\makeatother

\newcommand{\C}{\mathcal{C}}
\newcommand{\D}{\mathcal{D}}
\newcommand{\W}{\mathcal{W}}
\newcommand{\F}{\mathcal{F}}
\newcommand{\de}{\partial}
\newcommand{\p}{\times}
\renewcommand{\vec}{\overrightarrow}

\newtheorem*{thmN}{Theorem}

\newtheorem{thm}{Theorem}[section]
\newtheorem{prop}[thm]{Proposition}
\newtheorem{lem}[thm]{Lemma}
\newtheorem{cor}[thm]{Corollary}

\newtheorem{defnot}[thm]{Definition and notation}
\newcommand{\bdn}{\begin{defnot}}
	\newcommand{\edn}{\end{defnot}}
\newcommand{\bp}{\begin{prop}}
	\newcommand{\ep}{\end{prop}}
\newcommand{\bth}{\begin{thm}}
	\renewcommand{\eth}{\end{thm}}
\newcommand{\bpf}{\begin{proof}}
	\newcommand{\epf}{\end{proof}}
\newcommand{\bc}{\begin{cor}}
	\newcommand{\ec}{\end{cor}}

\theoremstyle{definition}
\newtheorem{defn}[thm]{Definition}
\newtheorem{rem}[thm]{Remark}
\newtheorem{exa}[thm]{Example}
\newcommand{\bd}{\begin{defn}}
\newcommand{\ed}{\end{defn}}
\newtheorem{nota}[thm]{Notation}

\newtheorem{qu}[thm]{Question}

\newcommand{\fL}[1]{\ar@{->}[ll]_-{#1}}
\newcommand{\fR}[1]{\ar@{->}[rr]^-{#1}}
\newcommand{\fRr}[1]{\ar@{->}[rrr]^-{#1}}
\newcommand{\fD}[1]{\ar@{->}[dd]_-{#1}}
\newcommand{\fU}[1]{\ar@{->}[uu]^-{#1}}
\newcommand{\f}[2]{\ar@{->}[#1]|{#2}}
\newcommand{\ff}[2]{\ar@2{->}[#1]|{#2}}
\newcommand{\frr}[1]{\ar@{->}[rrrr]^-{#1}}

\newcommand{\fl}[1]{\ar@{->}[l]_-{#1}}
\newcommand{\fr}[1]{\ar@{->}[r]^-{#1}}
\newcommand{\fd}[1]{\ar@{->}[d]_-{#1}}
\newcommand{\fu}[1]{\ar@{->}[u]^-{#1}}

\renewcommand{\top}{{\mathbf{Top}}}
\newcommand{\iso}{\cong}

\renewcommand{\leq}{\leqslant}
\renewcommand{\geq}{\geqslant}

\def\cartesien{%
	\ar@{-}[]+R+<6pt,-2pt>;[]+RD+<6pt,-6pt>%
	\ar@{-}[]+D+<2pt,-6pt>;[]+RD+<6pt,-6pt>%
}
\def\cocartesien{%
	\ar@{-}[]+L+<-6pt,+2pt>;[]+LU+<-6pt,+6pt>%
	\ar@{-}[]+U+<-2pt,+6pt>;[]+LU+<-6pt,+6pt>%
}
\def\hocartesien{%
	\ar@{-}[]+R+<6pt,-2pt>;[]+RD+<6pt,-6pt>_{h}%
	\ar@{-}[]+D+<2pt,-6pt>;[]+RD+<6pt,-6pt>%
}
\def\hococartesien{%
	\ar@{-}[]+L+<-6pt,+2pt>;[]+LU+<-6pt,+6pt>_{h}%
	\ar@{-}[]+U+<-2pt,+6pt>;[]+LU+<-6pt,+6pt>%
}

\newcommand{\brm}[1]{{\rm{\mathbf{#1}}}}
\newcommand{\dtop}{{\brm{Flow}}}
\newcommand{\set}{{\brm{Set}}}
\newcommand{\glob}{{\mathrm{Glob}}}
\newcommand{\ttt}{two-out-of-three property}
\DeclareMathOperator{\id}{Id}
\newcommand{\mins}[1]{\texorpdfstring{$#1$}{Lg}}
\newcommand{\liminj}{\varinjlim}

\renewcommand{\P}{\mathbb{P}}

\makeatletter
\def\varholim@#1#2{%
	\vtop{\m@th\ialign{##\cr
			\hfil$#1\operator@font holim$\hfil\cr
			\noalign{\nointerlineskip\kern1.5\ex@}#2\cr
			\noalign{\nointerlineskip\kern-\ex@}\cr}}%
}
\def\holimproj{%
	\mathop{\mathpalette\varholim@{\leftarrowfill@\textstyle}}\nmlimits@
}
\def\holiminj{%
	\mathop{\mathpalette\varholim@{\rightarrowfill@\textstyle}}\nmlimits@
}
\makeatother

\DeclareMathOperator{\cell}{{\brm{cell}}}
\DeclareMathOperator{\cof}{{\brm{cof}}}
\DeclareMathOperator{\inj}{{\brm{inj}}}
\newcommand{\ddownarrow}{{\downarrow}}
\newcommand{\Di}{\mathbf{D}}
\newcommand{\Sp}{\mathbf{S}}
\newcommand{\vd}{\vec{d}\!_1}
\DeclareMathOperator{\TT}{T}
\newcommand{\LMet}{\mathbf{LvMet}}
\newcommand{\SLMet}{\mathbf{PseudoMet}}

\allowdisplaybreaks

\begin{document}

\begin{abstract}
	Symmetric transverse sets were introduced to make the construction of the parallel product with synchronization for process algebras functorial. It is proved that one can do directed homotopy on symmetric transverse sets in the following sense. A q-realization functor from symmetric transverse sets to flows is introduced using a q-cofibrant replacement functor of flows. By topologizing the cotransverse maps, the cotransverse topological cube is constructed. It can be regarded both as a cotransverse topological space and as a cotransverse Lawvere metric space.   A natural realization functor from symmetric transverse sets to flows is introduced using Raussen's notion of natural $d$-path extended to symmetric transverse sets thanks to their structure of Lawvere metric space. It is proved that these two realization functors are homotopy equivalent on cofibrant symmetric transverse sets by using the fact that the small category defining symmetric transverse sets is c-Reedy in Shulman's sense. This generalizes to symmetric transverse sets results previously obtained for precubical sets. 
\end{abstract}
	
\maketitle
\setcounter{tocdepth}{1}
\tableofcontents
\hypersetup{linkcolor = dark-red}

\part*{Introduction}

\subsection*{Presentation}

Precubical sets are an important combinatorial model for directed homotopy \cite{DAT_book}. The $n$-cube represents the concurrent execution of $n$ actions. It has been known for a long time that the usual degeneracy maps used in non-directed homotopy theory are not convenient for directed homotopy. The purpose of this paper is to introduce a convenient notion of degeneracy map for doing directed homotopy. 

This is the second paper about \textit{symmetric transverse sets}. This notion is introduced in \cite{symcub} to make the construction of the parallel product with synchronization of process algebras functorial. It is proved in \cite[Theorem 3.1.15]{symcub} that it is the \textit{only} solution to achieve such a goal. A symmetric transverse set is a presheaf on the category $\widehat{\square}$ generated by the posets $[n]=\{0<1\}^n$ for $n\geq 0$ and by all \textit{cotransverse maps}. The latter are the strictly increasing maps preserving adjacency (Definition~\ref{cotransverse}). Note that to avoid cumbersome and inconsistent terminology, the words \textit{adjacency-preserving map} and \textit{transverse symmetric precubical set} of \cite{symcub} are replaced in this paper by \textit{cotransverse map} and \textit{symmetric transverse set} respectively. 

All coface maps and all symmetry maps are cotransverse. The latter cannot be regarded as degeneracy maps. However, a map like $\gamma_1:[2] \rightarrow [2]$ defined by $\gamma_1(\epsilon_1,\epsilon_2) = (\max(\epsilon_1,\epsilon_2),\min(\epsilon_1,\epsilon_2))$ is also cotransverse. Since $\gamma_1(1,0) = \gamma_1(0,1) = (1,0)$, the map $\gamma_1$ adds a degenerate $2$-cube by crushing the $2$-cube transversally to the direction of time which goes from $(0,0)$ to $(1,1)$. The \textit{transverse degeneracy maps} are exactly the cotransverse maps from $[n]$ to itself for $n\geq 2$ which are not one-to-one. All examples coming from computer science are symmetric transverse sets freely generated by precubical sets, and even most of them are freely generated by non-positively curved precubical sets \cite[Proposition~1.29]{zbMATH07226006}. Symmetric transverse sets provide a setting with degeneracies that fit with directed homotopy. They do not provide new examples for computer science: there does not seem to be any interpretation in computer science of a degenerate cube like the one given by $\gamma_1$.

The first goal of this paper is to prove that the notion of symmetric transverse set is a convenient framework for Raussen's notion of natural $d$-path. After translating the cotransverse maps into continuous maps thanks to a max-min formula, we obtain the cotransverse topological cube which is a cotransverse object both in the category of topological spaces and in the category of Lawvere metric spaces. The point is that the cotransverse maps are quasi-isometric: they preserve finite distances indeed. This implies that the topological version of the cotransverse maps takes natural $d$-paths of the topological cube to natural $d$-paths. This enables us to define a natural $d$-path between two vertices of a symmetric transverse set $K$ as a Moore composition of quasi-isometries from $\vec{[0,n]}$ for some integer $n\geq 1$ to the realization $|K|_{\vd}$ of $K$ as a Lawvere metric space. The following theorem summarizes the results of Part~\ref{metric}: 

\begin{thmN} (Theorem~\ref{cotransverse-cube}, Theorem~\ref{cotransverse-Lawverecube} and Corollary~\ref{samecube})
	For every cotransverse map $f:[m]\to [n]$, the map $\TT(f):[0,1]^m\to [0,1]^n$ of Definition~\ref{topologize1} gives rise to a quasi-isometry of Lawvere metric spaces for the Lawvere metric $\vd$ of Definition~\ref{d1} and enables us to define cotransverse objects both in the category of topological spaces and in the category of Lawvere metric spaces.
\end{thmN}

In \cite{NaturalRealization}, two realization functors from precubical sets to flows, a q-realization functor using a q-cofibrant replacement functor of the q-model structure of flows and the natural realization functor using Raussen's notion of natural $d$-path, are compared and are proved to be homotopy equivalent. The second goal of this paper is to generalize these results to symmetric transverse sets. The small category $\widehat{\square}$ is not Reedy. It is not Reedy in Berger-Moerdijk's sense \cite[Definition~1.1]{g-Reedy} or in Cisinski's sense \cite[Definition~8.1.1]{Cisinski-Book} either. However, it is c-Reedy in Shulman's sense \cite[Definition~8.25]{c-Reedy}. It is the key point to compare a q-realization and the natural realization of a symmetric transverse set as a flow. The difference with the setting of precubical sets studied in \cite{NaturalRealization} is that there is in general only a zigzag of natural transformations between the two realization functors and that the second part of the main theorem holds only for cofibrant symmetric transverse sets in the sense of Definition~\ref{cofibrant-transverse}. The following theorem summarizes the results of Part~\ref{homotopical}:  

\begin{thmN} (Theorem~\ref{exists}, Proposition~\ref{extend_nat}, Theorem~\ref{same-rea})
The natural realization functor $|-|_{nat}$ from symmetric transverse sets to flows defined in Definition~\ref{natural-rea} is an m-realization functor which extends the natural realization of precubical sets. Let $|-|_q$ be a q-realization functor of symmetric transverse sets. There exists an m-realization functor $F:\widehat{\square}^{op}\set \to \dtop$ and two natural transformations inducing bijections on the sets of states \[|-|_q\Longleftarrow F(-)\Longrightarrow |-|_{nat}\] such that for all cofibrant symmetric transverse sets $K$, and in particular for all symmetric transverse sets freely generated by a precubical set, and all $(\alpha,\beta)\in K_0\p K_0$, there is the zigzag of natural homotopy equivalences between m-cofibrant topological spaces 
\[\xymatrix
{
	\P_{\alpha,\beta} |K|_q & \fl{\simeq}  \P_{\alpha,\beta} F(K) \fr{\simeq} &\P_{\alpha,\beta} |K|_{nat}
}
.\] If $|-|_{q}$ is cofibrant as a q-realization functor, then one can suppose that $F= |-|_{q}$.
\end{thmN}

By Theorem~\ref{cofibrant-case}, the zigzag of homotopy equivalences between the spaces of execution paths would hold for non-cofibrant symmetric transverse sets (by considering a \textit{cofibrant} q-realization functor in the sense of Definition~\ref{homotopy-r-realization}) if the natural realization functor was cofibrant as well. It is unlikely that it is true but we cannot prove it. 

\subsection*{Future works} 

Subsequent papers will study the homotopical properties of transverse degeneracies. The terminology of \textit{symmetric} transverse set is used throughout the paper. The correct definition of a transverse set without symmetry maps is given in \cite{ThickCubes}. A convenient notion of \textit{labelled} symmetric transverse set should lead to a notion of higher dimensional transition system with degenerate higher dimensional transitions. The latter should yield model categories of higher dimensional transition systems which are not Quillen equivalent to discrete model categories, unlike e.g. in \cite[Theorem~7.5]{pastsim}. On the other hand, the \textit{non-symmetric} transverse sets defined in \cite{ThickCubes} should provide a convenient combinatorial setting for doing Directed Algebraic Topology.

\subsection*{Outline of the paper}

Part~\ref{metric} studies symmetric transverse sets from a metric point of view. 

Section~\ref{lawveremetric} recalls some basic facts about Lawvere metric spaces. The Lawvere metric $\vd$ (Definition~\ref{d1}) plays an important role in many places of the paper. 

Section~\ref{precubical-transverse} recalls some basic facts about precubical sets and symmetric transverse sets and the relations between one another. It also expounds in Proposition~\ref{free_square} a missing argument in the proof of \cite[Corollary~2.2.11]{symcub}. The section ends by introducing the notion of cofibrant symmetric transverse set in Definition~\ref{cofibrant-transverse} and by giving the basic properties. This notion has no analogue for precubical sets because, in some sense, all precubical sets are cofibrant. 

Section~\ref{topologize} starts from the observation made in Proposition~\ref{maxminxformula} to topologize the cotransverse maps and proves some useful properties about them. It culminates with Theorem~\ref{cotransverse-cube} and  \ref{cotransverse-Lawverecube} which expound the cotransverse topological cube and the cotransverse Lawvere cube.

Section~\ref{nat-path} is devoted to define the notion of natural $d$-path of a symmetric transverse set. It requires to recall what is the underlying topological space of a Lawvere metric space and to make some calculations about the cubes and more generally about the symmetric transverse sets. It is proved in Proposition~\ref{natural-metric} that natural $d$-paths of a topological cube are quasi-isometry for the Lawvere metric $\vd$. It enables us to define the natural $d$-paths of a symmetric transverse set as being locally, on each cube, a quasi-isometry. 

Part~\ref{homotopical} studies realization functors of symmetric transverse sets. 

Section~\ref{c-Reedy} proves that the category of cotransverse objects of a model category satisfying some mild conditions has a structure of a c-Reedy model category and that it coincides with the projective model structure. It enables us to give a necessary and sufficient condition for a cotransverse object to be projective cofibrant in Theorem~\ref{proj-cof-suff-cond}. The latter condition is used in Proposition~\ref{elm} to prove that the projective r-cofibrant replacement of the cotransverse flow associated with an r-realization functor gives rise to an r-realization functor. It is the key fact to prove Theorem~\ref{homotopy-natural} and Theorem~\ref{same-rea}. 

Section~\ref{rea_sec} defines the notions of q-realization, m-realization and h-realization of a symmetric transverse set as a flow and compares them both in the non-cofibrant case in Theorem~\ref{homotopy-natural} and in the cofibrant case in Theorem~\ref{cofibrant-case}. Theorem~\ref{exists} provides an example of a q-realization functor from symmetric transverse sets to flows. 

Section~\ref{concl} concludes this paper by defining the natural realization of a symmetric transverse set in Definition~\ref{natural-rea} and by comparing it in Theorem~\ref{same-rea} with a q-realization functor.

\subsection*{Prerequisites and notations}

The reading of \cite{NaturalRealization} is not required to understand this paper. A model category is a bicomplete category $\mathcal{M}$ equipped with a class of cofibrations $\C$, a class of fibrations $\F$ and a class of weak equivalences $\W$ such that: 1) $\W$ is closed under retract and satisfies \ttt, 2) the pairs $(\C,\W\cap \F)$ and $(\C\cap \W,\F)$ are functorial weak factorization systems. We refer to \cite[Chapter~1]{MR99h:55031} and to \cite[Chapter~7]{ref_model2} for the basic notions about general model categories. We refer to \cite{TheBook} for locally presentable categories and to \cite{MR2506258} for combinatorial model categories.

Let $\mathcal{I}$ be a set of maps of a cocomplete category $\C$. The notation $f\boxslash g$ means that $g$ satisfies the right lifting property (RLP) with respect to $f$; ${^\boxslash}\mathcal{I} = \{g, \forall f \in \mathcal{I}, g\boxslash f\}$; $\mathcal{I}^\boxslash = \inj(\mathcal{I}) = \{g, \forall f \in \mathcal{I}, f\boxslash g\}$; $\cof(\mathcal{I})={^\boxslash}(\mathcal{I}^\boxslash)$; $\cell(\mathcal{I})$ is the class of transfinite compositions of pushouts of elements of $\mathcal{I}$. A \textit{cellular} object (with respect to $\mathcal{I}$) is an object $X$ such that the canonical map $\varnothing \to X$ belongs to $\cell(\mathcal{I})$. A \textit{cofibrant} object (with respect to $\mathcal{I}$) is an object $X$ such that the canonical map $\varnothing \to X$ belongs to $\cof(\mathcal{I})$. By \cite[Corollary~2.1.15]{MR99h:55031}, if the category $\C$ is locally presentable, then the cofibrant objects with respect to $\mathcal{I}$ are exactly the retracts of the cellular objects with respect to $\mathcal{I}$. 

The set of maps from $X$ to $Y$ of a category $\C$ is denoted by $\C(X,Y)$. $\varnothing$ denotes the initial object and $\mathbf{1}$ the final object of a category. $\set$ is the category of sets with all set maps. $\C^I$ is the category of functor from a small category $I$ to a category $\C$ together with the natural transformations. $\iso$ means isomorphism, $\simeq$ means weak equivalence or homotopy equivalence, depending on the context. For an object $X$ of a category $\C$ and a set $S$, $S.X$ denotes the sum $\coprod_S X$ and $X^S$ denotes the product $\prod_S X$. 

The category $\top$ denotes either the category of \textit{$\Delta$-generated spaces} or of \textit{$\Delta$-Hausdorff $\Delta$-generated spaces} (cf. \cite[Section~2 and Appendix~B]{leftproperflow}).

\subsection*{Warnings}

All $d$-paths considered in this paper are \textit{tame} in the sense of \cite[Section~2.9]{MR4070250} and \textit{nonconstant}. These adjectives are therefore understood everywhere. See also the end of Section~\ref{nat-path} where the notions of $d$-path and natural $d$-path of a symmetric transverse set are introduced.

\subsection*{Acknowledgments}

I am very grateful to the anonymous referee for the extremely detailed report.

\part{Metric study of symmetric transverse sets}
\label{metric}

\section{The Lawvere directed \mins{n}-cube}
\label{lawveremetric}

Since there are several variants of the notion of metric space in the mathematical literature, the one which is used in this paper is recalled. The symmetric version will have to be recalled in Section~\ref{nat-path}.

\bd \cite{LawvereMetric}
A \textit{Lawvere metric space} $(X,d)$ is a set $X$ equipped with a map $d:X\p X\to [0,\infty]$ called a \textit{(Lawvere) metric} such that:
\begin{itemize}
	\item $\forall x\in X,d(x,x)=0$
	\item $\forall (x,y,z)\in X\p X\p X, d(x,y)\leq d(x,z)+d(z,y)$.
\end{itemize}
A map $f:(X,d)\to (Y,d)$ of Lawvere metric spaces is a set map $f:X\to Y$ which is \textit{non-expansive}, i.e. $\forall (x,y)\in X\p X, d(f(x),f(y))\leq d(x,y)$. A non-expansive map $f:(X,d)\to (Y,d)$ is \textit{quasi-isometric} if $\forall (x,y)\in X\p X, d(x,y)<\infty \Rightarrow d(f(x),f(y))= d(x,y)$.
\ed

\begin{nota}
	The category of Lawvere metric spaces is denoted by $\LMet$.
\end{nota}

The category of Lawvere metric spaces is bicomplete since it is the category of small categories enriched over $([0,\infty],\geq,+,0)$ \cite{LawvereMetric}.

\begin{nota}
	Let $[0] = \{()\}$ and $[n] = \{0,1\}^n$ for $n \geq 1$. By convention, one has $\{0,1\}^0=[0]=\{()\}$. In the sequel, for all $n\geq 1$, both the sets $[n]$ and $[0,1]^n$ are equipped with the product order. By convention, $[0,1]^0$ is a singleton.
\end{nota}

\bd \label{d1}
	Let $x=(x_1,\dots,x_n)$ and $x'=(x'_1,\dots,x'_n)$ be two elements of $[0,1]^n$ with $n\geq 1$. Let $\vd:[0,1]^n\p [0,1]^n \to [0,\infty]$ be the set map defined by
	\[
	\vd(x,x') = 
	\begin{cases}
		\displaystyle\sum\limits_{i=1}^{n} |x_i-x'_i| & \hbox{ if } x\leq x'\\
		\infty & \hbox {otherwise.}
	\end{cases}
	\]
	For $n=1$, it is $\vec{[0,1]}$ of \cite[Example~3.2]{zbMATH07226006}.
\ed

\bp
Let $n\geq 0$. The set map $\vd:[0,1]^n\p [0,1]^n\to [0,\infty]$ is a Lawvere metric. It restricts to a Lawvere metric on $\{0,1\}^n$.
\ep

\bpf
Let $x,y,z\in [0,1]^n$. If $\vd(x,z)+\vd(z,y)$ is finite, then $x\leq z\leq y$, which implies that $\vd(x,y)$ is finite and that $\vd(x,y)=\vd(x,z)+\vd(z,y)$. If $\vd(x,z)+\vd(z,y)$ is infinite, then the inequality $\vd(x,y)\leq \vd(x,z)+\vd(z,y)$ always holds.
\epf

\section{Precubical and symmetric transverse set}
\label{precubical-transverse}

\begin{nota}
	Let $A\subset \{1,\dots,n\}$. Denote by $\epsilon_A$ the tuple $(\epsilon_1,\dots,\epsilon_n)$ with $\epsilon_i=0$ if $i\notin A$ and $\epsilon_i=1$ of $i\in A$. Let $0_n=\epsilon_\varnothing$ and $1_n=\epsilon_{\{1,\dots,n\}}$.
\end{nota}

Let $\delta_i^\alpha : [n-1] \rightarrow [n]$ be the \textit{coface map} defined for $1\leq i\leq n$ and $\alpha \in \{0,1\}$ by $\delta_i^\alpha(\epsilon_1, \dots, \epsilon_{n-1}) = (\epsilon_1,\dots, \epsilon_{i-1}, \alpha, \epsilon_i, \dots, \epsilon_{n-1})$. The small category $\square$ is by definition the subcategory of the category of posets with the set of objects $\{[n],n\geq 0\}$ and generated by the morphisms $\delta_i^\alpha$. The maps of $\square$ are called the \textit{cocubical} maps. 

\bd \cite{Brown_cube} The category of presheaves over $\square$, denoted by $\square^{op}\set$, is called the category of \textit{precubical sets}. Let $\square[n]:=\square(-,[n])$.  For $K\in \square^{op}\set$, denote by  $K_n=K([n])$ the set of \textit{$n$-cubes} of $K$.  For $c\in K_n$, let $n=\dim(c)$. Let $f:[m]\to [n]$ be a cocubical map. It gives rise to a set map denoted by $f^*:K_n\to K_m$. An element of $K_0$ is called a \textit{vertex} of $K$. \ed

The following definition is equivalent to \cite[Definition~2.1.5]{symcub}. 

\bd \label{cotransverse} A set map $f:[m] \rightarrow [n]$ is \textit{cotransverse} if it is strictly increasing and if $\forall x,y\in [m], \vd(x,y)=1 \hbox{ implies }\vd(f(x),f(y)) = 1$.
\ed

The adjective adjacency-preserving is used in \cite{symcub} instead. The word cotransverse is preferred because it is consistent with the terminology of symmetric transverse sets~\footnote{Unlike in \cite{symcub}, the word precubical is omitted.}.

By \cite[Proposition~2.1.6]{symcub}, for any $n\geq 1$, the coface map $\delta_i^\alpha : [n-1] \rightarrow [n]$ is cotransverse and any strictly increasing map from $[n]$ to itself is cotransverse as well. Let $\sigma_i:[n] \rightarrow [n]$ be the set map defined for $1\leq i\leq n-1$ and $n\geq 2$ by $\sigma_i(\epsilon_1, \dots, \epsilon_{n}) = (\epsilon_1, \dots, \epsilon_{i-1},\epsilon_{i+1},\epsilon_{i}, \epsilon_{i+2},\dots,\epsilon_{n})$. These maps are called the \textit{symmetry maps} \cite{MR1988396}. The symmetry maps are clearly cotransverse.

\begin{nota}
	Let $\widehat{\square}$ be the small subcategory of the category of posets generated by the cotransverse maps.
\end{nota}

The following proposition is crucial in many places of this paper.

\bp \cite[Proposition~3.1.14]{symcub} \label{decomposition_distance} Let $0\leq m\leq n$. Every cotransverse (resp. cotransverse one-to-one) map $f:[m] \rightarrow [n]$ factors uniquely as a composite $[m] \stackrel{\psi}\longrightarrow [m] \stackrel{\phi}\longrightarrow [n]$ with $\phi\in \square$ and $\psi$ cotransverse (resp. cotransverse one-to-one).  \ep

By a cardinality argument, if $\psi:[m] \rightarrow [m]$ is one-to-one, then it is bijective and therefore it is a symmetry map. Thus the one-to-one cotransverse maps are composites of coface maps and symmetry maps in a unique way. 

\bd \cite[Definition~2.1.13]{symcub} The category of presheaves over $\widehat{\square}$, denoted by $\widehat{\square}^{op}\set$, is called the category of \textit{symmetric transverse sets}. Let $\widehat{\square}[n]:=\widehat{\square}(-,[n])$. For $K\in \widehat{\square}^{op}\set$, denote by  $K_n=K([n])$ the set of \textit{$n$-cubes} of $K$.  For $c\in K_n$, let $n=\dim(c)$. Let $f:[m]\to [n]$ be a cotransverse map. It gives rise to a set map denote by $f^*:K_n\to K_m$. An element of $K_0$ is called a \textit{vertex} of $K$.
\ed

\bd Let $\C$ be a category. A \textit{cotransverse object} of $\C$ is a functor $\widehat{\square}\to \C$.
\ed

There is the elementary proposition: 

\bp \label{eq-cat}
Let $\C$ be a cocomplete category. Let $X:\widehat{\square}\to \C$ be a cotransverse object of $\C$. Let \[\widehat{X}(K)=\int^{[n]\in\widehat{\square}} K_n.X([n])\] The mapping $X\mapsto \widehat{X}$ induces an equivalence of categories between the category of cotransverse objects of $\C$ and the colimit-preserving preserving functors from $\widehat{\square}^{op}\set$ to $\C$.
\ep

Proposition~\ref{eq-cat} is a particular case of \cite[Remark~3.2.3]{coend-calculus}, a cotransverse object in a cocomplete category being an example of a \textit{nerve-realization context} as defined in \cite[Definition~3.2.1]{coend-calculus}. Proposition~\ref{eq-cat} holds for any small category instead of $\widehat{\square}$. For the case of cosimplicial objects, see \cite[Proposition~3.1.5]{MR99h:55031}. The reader might also be interested in reading \cite[Remark~6.5.9]{CategoryInContext} for another presentation of the general case.

\begin{nota}
	For the sequel, the cotransverse object associated with a colimit-preserving functor $F:\widehat{\square}^{op}\set\to \C$ is denoted by $F(\widehat{\square}[*])$. 
\end{nota}

\begin{nota} \label{def_h}
	Let $n\geq 1$. Let $h:[0,1]^n \to [0,n]$ be the continuous map defined by \[h(x_1,\dots,x_n) = \sum_{i=1}^n x_i.\]
	Note that for all $x,y\in [0,1]^n$, $x\leq y$ implies $h(x)\leq h(y)$ and that $x\leq y$ and $h(x)=h(y)$ implies $x=y$.
\end{nota}

\bp \label{preserv-h-dis}
Let $n\geq 1$. Let $f:[n]\to [n]$ be a cotransverse map. Then for all $(\epsilon_1,\dots,\epsilon_n)\in [n]$, one has $h(\epsilon_1,\dots,\epsilon_n) = h(f(\epsilon_1,\dots,\epsilon_n))$.
\ep

\bpf
We proceed by induction on $h(\epsilon_1,\dots,\epsilon_n)$. Consider the increasing sequence \[\epsilon_\varnothing< \epsilon_{\{1\}}<\epsilon_{\{1,2\}}<\dots< \epsilon_{\{1,2,\dots,n\}}\] of elements of $[n]$. The map $f$ being cotransverse by hypothesis, one has 
\begin{align*}
	&\vd\big(f(\epsilon_\varnothing),f(\epsilon_{\{1\}})\big) = 1,\\
	&\vd\big(f(\epsilon_{\{1\}}),f(\epsilon_{\{1,2\}})\big) = 1,\\
	& \dots \\
	&\vd\big(f(\epsilon_{\{1,\dots,n-1\}}),f(\epsilon_{\{1,\dots,n\}})\big) = 1.
\end{align*}
Since $f:[n]\to [n]$ is strictly increasing, we obtain $f(\epsilon_\varnothing)=\epsilon_\varnothing$ and $f(\epsilon_{\{1,\dots,n\}}) = \epsilon_{\{1,\dots,n\}}$. We deduce that $h(\epsilon_\varnothing)=hf(\epsilon_\varnothing)$ and that $h(\epsilon_{\{1,\dots,n\}})=hf(\epsilon_{\{1,\dots,n\}})$. The formula is therefore proved for $h(\epsilon_1,\dots,\epsilon_n)=0$ (and also for $h(\epsilon_1,\dots,\epsilon_n)=n$). Suppose the formula proved for all $(\epsilon_1,\dots,\epsilon_n)\in [n]$ such that $h(\epsilon_1,\dots,\epsilon_n)\leq H<n$. Let $(\epsilon_1,\dots,\epsilon_n)\in [n]$ such that $h(\epsilon_1,\dots,\epsilon_n)=H+1 \geq 1$. There exists $(\epsilon'_1,\dots,\epsilon'_n)\in [n]$ with $h(\epsilon'_1,\dots,\epsilon'_n)=H$ and $(\epsilon'_1,\dots,\epsilon'_n)<(\epsilon_1,\dots,\epsilon_n)$. We deduce that $\vd((\epsilon'_1,\dots,\epsilon'_n),(\epsilon_1,\dots,\epsilon_n))=1$. The map $f$ being cotransverse, we obtain $\vd(f(\epsilon'_1,\dots,\epsilon'_n),f(\epsilon_1,\dots,\epsilon_n))=1$. We obtain the equalities $h(f(\epsilon_1,\dots,\epsilon_n))=h(f(\epsilon'_1,\dots,\epsilon'_n))+1 = H+1$, the first equality by definition of $\vd$ and the second equality by induction hypothesis. 
\epf

As a corollary, we obtain the following proposition.

\bp \label{pseudoiso-dis}
Let $\psi:[m]\to [n]$ be a cotransverse map. Then $\psi$ induces a map of Lawvere metric spaces from $[m]$ to $[n]$ which is quasi-isometric.
\ep

A cotransverse map is not necessarily an isometry. For example, the map $\gamma_1:[2] \rightarrow [2]$ defined by $\gamma_1(\epsilon_1,\epsilon_2) = (\max(\epsilon_1,\epsilon_2),\min(\epsilon_1,\epsilon_2))$ is cotransverse and $\gamma_1(1,0) = \gamma_1(0,1) = (1,0)$. Note that $\vd((0,1),(1,0))=\infty$.

\begin{nota} \label{free-transverse}
	The inclusion of small categories $\square\subset \widehat{\square}$ induces a forgetful functor $\widehat{\omega}:\widehat{\square}^{op}\set \to \square^{op}\set$ which has a left adjoint $\widehat{\mathcal{L}}:\square^{op}\set \to \widehat{\square}^{op}\set$ which is called the \textit{free symmetric transverse set} generated by a precubical set. 
\end{nota}

\bp \label{restriction}
	For a precubical (symmetric transverse resp.) set $K$, the data
	\[
	(K_{\leq n})_p = \begin{cases}
		K_p & \hbox{ if }p\leq n\\
		\varnothing & \hbox{ if }p> n.
	\end{cases}
	\]
	assemble into a precubical (symmetric transverse resp.) set denoted by $K_{\leq n}$. Moreover, the functor $K \mapsto K_{\leq n}$ is colimit-preserving.
\ep 

\bpf
The first part is due to the fact that $\square([m],[n])=\widehat{\square}([m],[n]) = \varnothing$ when $m>n$. The second part is due to the fact that colimits of presheaves are calculated objectwise.
\epf

\begin{nota}
	Let $\de\square[n] = \square[n]_{\leq n-1}$ and $\de\widehat{\square}[n]=\widehat{\square}[n]_{\leq n-1}$ for all $n\geq 0$.
\end{nota}

\begin{figure}
	\[{
		\xymatrix@C=3em@R=3em{ [p]
			\ar@{=}[rd] \ar@{->}[dd]_-{k} \ar@{->}[rr]^-{\psi} && [p]
			\ar@{->}[rd]^-{h'g} \ar@{->}'[d][dd]^-{\psi'}  &\\
			& [p] \ar@{->}[rr]^(0.3){\psi_1} \ar@{->}[dd]_(0.3){k} && [r] \ar@{->}[dd]^-{\psi_3} \\
			[n] \ar@{=}[rd] \ar@{=}'[r][rr] &&
			[n] \ar@{-->}[rd]^-{\psi_2}& \\
			& [n] \ar@{->}[rr]^(0.3){\psi_2} &&
			[s]
		}
	} \]
	\caption{$(k\ddownarrow L)$ is connected}
	\label{terminal}
\end{figure}
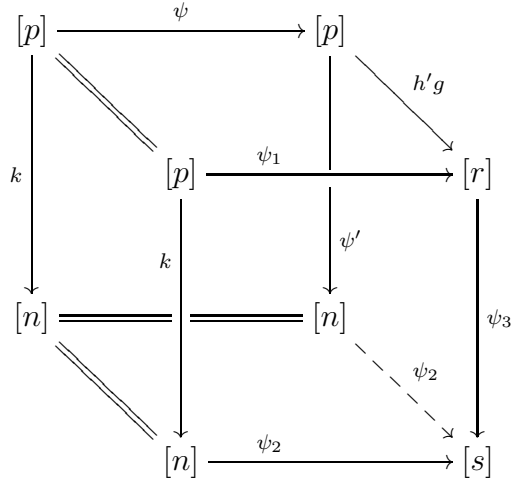

\bp \label{free_square}
For all $n\geq 0$, one has the isomorphism of symmetric transverse sets \[\widehat{\mathcal{L}}(\square[n]) \iso \widehat{\square}[n].\] There is the isomorphism of symmetric transverse sets \[\widehat{\mathcal{L}}(\de\square[n]) \iso \de\widehat{\square}[n]\] for all $n\geq 0$. 
\ep

\bpf
The first statement is \cite[Proposition 2.1.14]{symcub}. The short argument is repeated for the ease of the reader. For every symmetric transverse set $K$, one has $K_n=(\widehat{\omega}K)_n$ for all $n\geq 0$. Since the functor $\square\subset \widehat{\square}$ is the identity on objects, we obtain for all $n\geq 0$ the bijections 
\[\widehat{\square}^{op}\set(\widehat{\mathcal{L}}(\square[n]),K) \iso \square^{op}\set(\square[n],\widehat{\omega} K) = (\widehat{\omega} K)_n = K_n = \widehat{\square}^{op}\set(\widehat{\square}[n],K).\] By Yoneda's lemma, one
obtains the isomorphism $\widehat{\mathcal{L}}(\square[n]) \iso \widehat{\square}[n]$ for all $n\geq 0$. The second statement is stated with an incorrect argument in the proof of \cite[Corollary~2.2.11]{symcub}. The missing argument is explained now. Consider the small category $J'$ such that the objects are the coface maps $[p]\to [n] \in \square$ with $p<n$ and such that the morphisms of $J'$ are the commutative squares of the form 
\[
\xymatrix@C=4em@R=4em
{
	[p] \fr{} \fd{\in \square} & [n] \ar@{=}[d] \\
	[q] \fr{}  & [n]
}
\]
Since $\square([p],[n]) = \varnothing$ for $p>n$ and since $\widehat{\mathcal{L}}:\square^{op}\set\to \widehat{\square}^{op}\set$ is colimit-preserving, we obtain the isomorphism of symmetric transverse sets \[\liminj_{[p]\to [n] \in J'} \widehat{\square}[p] \iso \widehat{\mathcal{L}}(\de\square[n]).\]
Consider the small category $J$ such that the objects are the maps $[p]\to [n] \in \widehat{\square}$ with $p<n$ and such that the morphisms are the commutative squares of the form 
\[
\xymatrix@C=4em@R=4em
{
	[p] \fr{} \fd{\in \widehat{\square}} & [n] \ar@{=}[d] \\
	[q] \fr{}  & [n]
}
\]
Since $\widehat{\square}([p],[n]) = \varnothing$ for $p>n$, we obtain the isomorphism of symmetric transverse sets \[\liminj_{[p]\to [n] \in J} \widehat{\square}[p] \iso \de\widehat{\square}[n].\]
Consider the inclusion functor $L:J'\to J$. It induces a map of symmetric transverse sets \[\widehat{\mathcal{L}}(\de\square[n]) \longrightarrow \de\widehat{\square}[n].\] By \cite[Theorem~1 p.~213]{MR1712872}, it suffices to prove that the comma category $(k\ddownarrow L)$ is nonempty and connected for all objects $k$ of $J$ to complete the proof. Let $k:[p]\to [n]$ be an object of $J$. We see immediately that the comma category $(k\ddownarrow L)$ is nonempty because it contains the commutative square
\[
\xymatrix@C=4em@R=4em
{
	[p] \fr{\psi} \fd{k} & [p] \ar@{->}[d]^-{\psi'\in\square} \\
	[n] \ar@{=}[r]  & [n]
}
\]
where the top map $\psi:[p]\to [p]$ is given by the unique factorization given by Proposition~\ref{decomposition_distance} of $k:[p]\to [n]$ as the composite of a map of $\widehat{\square}([p],[p])$ followed by a coface map $\psi'$. Consider another object 
\[
\xymatrix@C=4em@R=4em
{
	[p] \fr{\psi_1} \fd{k} & [r] \ar@{->}[d]^-{\psi_3\in\square} \\
	[n] \fr{\psi_2} & [s]
}
\]
of the comma category $(k\ddownarrow L)$. Consider the following diagram of $\widehat{\square}$:
\[
\xymatrix@C=5em@R=4em
{
[p] \ar@{=}[d]\fr{\psi}&[p] \fd{g}\fr{\psi'\in \square}& [n]\fd{\psi_2}\\
	[p] \ar@{=}[d]\ar@{->}@/^10pt/[r]^-{g\psi} \ar@{->}@/_10pt/[r]^-{h}&[p] \fd{h'\in \square}\fr{g' \in \square}& [s] \ar@{=}[d]\\
	[p] \fr{\psi_1} & [r] \fr{\psi_3\in \square} & [s]
}
\]
where the factorizations $\psi_2\psi'=g'g$ and $\psi_1=h'h$ are given by the factorization of Proposition~\ref{decomposition_distance}. We obtain $(\psi_3h')h=\psi_3\psi_1=\psi_2\psi'\psi=g'(g\psi)$. By uniqueness of the factorization of Proposition~\ref{decomposition_distance}, we deduce that $\psi_3h'=g'$ and $h=g\psi$. We deduce the map of $(k\ddownarrow L)$ depicted in Figure~\ref{terminal}. We conclude that the comma category $(k\ddownarrow L)$ is connected.
\epf

\begin{rem}
	In fact, we could prove that the comma category $(k\ddownarrow L)$ has an initial object given by the factorization of $k$ using Proposition~\ref{decomposition_distance}.
\end{rem}

\bd \label{cofibrant-transverse}
A symmetric transverse set $K$ is \textit{cellular} if the canonical map $\varnothing\to K$ is a transfinite composition of pushouts of the maps $\de\widehat{\square}[n]\to \widehat{\square}[n]$ for $n\geq 0$. Note that the map $\de\widehat{\square}[0]\to \widehat{\square}[0]$ is the map $C:\varnothing\to \{0\}$. A symmetric transverse set $K$ is \textit{cofibrant} if it is a retract of a cellular symmetric transverse set.
\ed

The category of symmetric transverse sets is locally presentable by \cite[Corollary~1.54]{TheBook}, being a presheaf category. Besides, the symmetric transverse sets $\de\widehat{\square}[n]$ are cofibrant for all $n\geq 0$ by Proposition~\ref{cp_cofibrant}. In other terms, the two sets of maps $\{\de\widehat{\square}[n]\to \widehat{\square}[n]\mid n\geq 0\}$ and $\{\de\widehat{\square}[n]\to \widehat{\square}[n]\mid n\geq 0\} \cup \{R:\{0,1\}\to \{0\}\}$ are tractable. By \cite[Theorem~1.4]{henry2020minimal}, the cofibrant symmetric transverse sets are therefore the cofibrant objects of the minimal model categories generated by these two sets of maps. It is not clear at this point whether $R:\{0,1\}\to \{0\}$ must be added or not to the set of generating cofibrations to have a non-trivial model category on the category of symmetric transverse sets. Besides, \cite[Theorem~1.4]{henry2020minimal} does not provide any geometric information. It is known by \cite[Corollary~4.10]{Nonunital} that removing $R:\{0,1\}\to \{0\}$ from the generating cofibrations of the q-model structure of flows (see Definition~\ref{def-flow} and Theorem~\ref{three}) leads to a minimal category without homotopy on the category of flows. However, the $(n+1)$-dimension globe $\glob(\Di^n)$ has two distinguished states whereas the $(n+1)$-cube has $2^{n+1}$ states. Thus, the induction which leads to \cite[Corollary~4.10]{Nonunital} does not work in the transverse case.

The terminal symmetric transverse set $T$ is not cofibrant. Indeed, $T_n$ is a singleton $\{c_n\}$ for all $n\geq 0$. If it was a retract of a cellular symmetric transverse set $K$, then the identity of $T$ would factor as a composite $T\to K\to T$. For all cotransverse maps $f:[n]\to [n]$, $f^*(c_n)=c_n$. However $f^*$ has no fix point in $K$: contradiction.

\bp \label{cp_cofibrant}
Let $K$ be a precubical set. Then the symmetric transverse set $\widehat{\mathcal{L}}(K)$ freely generated by $K$ is cellular. In particular, for all $n\geq 0$, the symmetric transverse sets $\de\widehat{\square}[n]$ and $\widehat{\square}[n]$ are cellular for all $n\geq 0$. 
\ep

\bpf
Let $K$ be a precubical set. Construct a transfinite tower $(K^\alpha)_{\alpha\geq 0}$ with $K^0=\varnothing$ and such that there is a map of precubical sets $K^\alpha\to K$ one-to-one on cubes as follows. Suppose that $K^\alpha$ is constructed. If there exists $n\geq 0$ such that $K^\alpha_n\to K_n$ is not onto, then take the least $n$ satisfying this property. Then there exists a pushout $K^\alpha \to K^{\alpha+1}$ of $\square[n]_{\leq n-1} \to \square[n]$ such that $K^{\alpha+1}\to K$ is still one-to-one on cubes because $\square([n],[n])$ is a singleton for all $n\geq 0$ and because colimits are calculated objectwise on presheaves. The transfinite induction stops eventually for a cardinality reason. Thus the proposition is a consequence of Proposition~\ref{free_square} and of the fact that the functor $\widehat{\mathcal{L}}:\square^{op}\set \to \widehat{\square}^{op}\set$ is colimit-preserving, being a left adjoint. 
\epf

\begin{figure}
	{\[
		\xymatrix{
			& (1,0,0) \ar@{->}[r] \ar@{->}[rrd] & (1,1,0) \ar@{->}[drr] && \\
			(0,0,0) \ar@{->}[r] \ar@{->}[ur] \ar@{->}[dr] & (0,1,0) \ar@{->}[ur]
			\ar@{->}[dr] && (1,0,1) \ar@{->}[r] &
			(1,1,1) \\
			& (0,0,1) \ar@{->}[r] \ar@{->}[rru] & (0,1,1) \ar@{->}[urr] && \\
			& (0,0,1) \ar@{->}[r] \ar@{->}[rrd] & (0,1,1) \ar@{->}[drr] && \\
			(0,0,0) \ar@{->}[r] \ar@{->}[ur] \ar@{->}[dr] & (0,0,1) \ar@{->}[ur]
			\ar@{->}[dr] && (1,0,1) \ar@{->}[r] &
			(1,1,1). \\
			& (0,0,1) \ar@{->}[r] \ar@{->}[rru] & (0,1,1) \ar@{->}[urr] && 
		}
		\] }
	\caption{The cotransverse map $f:[3]\to [3]$}
	\label{source}
\end{figure}

\bp
There exists a cofibrant symmetric transverse set which is not freely generated by a precubical set. 
\ep

\bpf
Consider the cotransverse map $f:[3]\to [3]$ defined as follows (it is the example \cite[Figure~5]{symcub} and it is depicted in Figure~\ref{source}): 
\[
f(x,y,z) = \begin{cases}
	(x,y,z) & \hbox{ for }(x,y,z)\in \{(0,0,0),(1,0,1),(1,1,1)\}\\
	(0,1,1) & \hbox{ for }(x,y,z)\in \{(1,1,0),(0,1,1)\\
	(0,0,1) & \hbox{ for } h(x,y,z) = x+y+z = 1
\end{cases}
\]
We immediately see that $f$ is cotransverse, each arrow adding exactly $1$ to the sum of elements of the triple. The map $f$ induces a map of symmetric transverse sets $\de f:\de\widehat{\square}[3]\to \de\widehat{\square}[3]$ which is not the image by $\widehat{\mathcal{L}}$ of a map of precubical sets from $\de{\square}[3]$ to $\de{\square}[3]$ because e.g. the $2$-dimensional subcube $(*,*,0)$ is crushed by $\de f$ to the concatenation of two edges $(0,0,0)\to (0,0,1) \to (0,1,1)$. Consider the pushout diagram of symmetric transverse sets 
\[
\xymatrix@C=3em@R=3em
{\de\widehat{\square}[3] \fd{\subset}\fr{\de f} & \de\widehat{\square}[3] \fd{} \\
	\widehat{\square}[3] \fr{} & \cocartesien X}
\]
Then the symmetric transverse set $X$ is cofibrant and it is not freely generated by a precubical set because it contains a degenerate $2$-cube.
\epf

\bp \label{finite_induction_cofibrant}
Let $K$ be a symmetric transverse set. It is cellular if and only if for all $n\geq 0$, there is the pushout diagram of symmetric transverse sets 
\[
\xymatrix@C=5em@R=3em
{
	\displaystyle\coprod\limits_{x\in \cell_{n+1}(K)} \widehat{\square}[n+1]_{\leq n} \fd{}\fr{} & K_{\leq n}\fd{} \\
	\displaystyle\coprod\limits_{x\in \cell_{n+1}(K)} \widehat{\square}[n+1] \fr{}  & \cocartesien K_{\leq n+1}
}
\]
where $\cell_{n}(K)$ is the set of $n$-cubes in the cellular decomposition of $K$.
\ep 

\bpf
If $K$ satisfies the property of the proposition, and since $K=\liminj K_{\leq n}$, then $K$ is cellular. Conversely, suppose that $K$ is a cellular symmetric transverse set. By Proposition~\ref{restriction}, the restriction functor $K\mapsto K_{\leq n}$ is colimit-preserving for all $n\geq 0$. Thus each natural map $\varnothing\to K_{\leq n}$ is a transfinite composition of pushouts of the maps $\de\widehat{\square}[p]_{\leq n}\to \widehat{\square}[p]_{\leq n}$ for $p\geq 0$. The point is that for all $p>n$, the map $\de\widehat{\square}[p]_{\leq n}\to \widehat{\square}[p]_{\leq n}$ is the identity of $\widehat{\square}[p]_{\leq n}$ by definition of $\de\widehat{\square}[p]_{\leq n}$. Moreover, each map of symmetric transverse sets $\widehat{\square}[n+1]_{\leq n} \to K$ factors uniquely as a composite $\widehat{\square}[n+1]_{\leq n} \to K_{\leq n}\to K$. Hence the proof is complete. 
\epf

\section{Cotransverse topological cube}
\label{topologize}

The purpose of this section is to topologize the cotransverse maps, more precisely to extend any cotransverse map $f:[m]\to [n]$ to a map of Lawvere metric spaces $\TT(f)$ from $([0,1]^m,\vd)$ to $([0,1]^n,\vd)$ which is quasi-isometric. The starting point is the following observation.

\bp \label{maxminxformula}
Let $n\geq 1$. Let $f=(f_1,\dots,f_n):[n]\to [n]$ be a cotransverse map. Then there is the equality 
\[
f_i(x_1,\dots,x_n) = \max_{(\epsilon_1,\dots,\epsilon_n)\in f_i^{-1}(1)} \min \{x_k\mid \epsilon_k=1\}
\]
for all $1\leq i\leq n$.
\ep

\bpf
There are two mutually exclusive cases: $f_i(x_1,\dots,x_n)=0$ or $f_i(x_1,\dots,x_n)=1$. Let us treat the case $f_i(x_1,\dots,x_n)=0$ at first. For all $(\epsilon_1,\dots,\epsilon_n)\in f_i^{-1}(1)$, $\min \{x_k\mid \epsilon_k=1\} = 1$ implies $(x_1,\dots,x_n)\geq (\epsilon_1,\dots,\epsilon_n)$, which implies $f_i(x_1,\dots,x_n)=1$: contradiction. Thus $f_i(x_1,\dots,x_n)=0$ implies that for all $(\epsilon_1,\dots,\epsilon_n)\in f_i^{-1}(1)$, one has $\min \{x_k\mid \epsilon_k=1\} = 0$. Assume now that $f_i(x_1,\dots,x_n)=1$. Then $(x_1,\dots,x_n)\in f_i^{-1}(1)$. Since $\min \{x_i\mid x_i=1\} = 1$, the proof is complete. 
\epf

To give the reader the intuition of Proposition~\ref{maxminxformula}, consider the cotransverse map $f:[3]\to [3]$ described in Figure~\ref{source}. Let $f=(f_1,f_2,f_3)$. The reader must keep in mind that, for boolean values, there are the equalities 
\[
\min(x,y) = x\hbox{ and }y \hbox{, }\max(x,y) = x\hbox{ or }y.
\]
If $x_1=1$ and $x_3=1$, or $x_1=1$ and $x_2=1$ and $x_3=1$, then $f_1(x_1,x_2,x_3)=1$. Thus \[f_1(x_1,x_2,x_3)=\max(\min(x_1,x_3),\min(x_1,x_2,x_3)).\] If $x_1=1$ and $x_2=1$, or $x_2=1$ and $x_3=1$, or $x_1=1$ and $x_2=1$ and $x_3=1$, then $f_2(x_1,x_2,x_3)=1$. Thus \[f_2(x_1,x_2,x_3)=\max(\min(x_1,x_2),\min(x_2,x_3),\min(x_1,x_2,x_3)).\] Finally, if $x_1=1$ and $x_2=1$, or $x_1=1$ and $x_3=1$, or $x_2=1$ and $x_3=1$, or $x_1=1$ and $x_2=1$ and $x_3=1$, then $f_3(x_1,x_2,x_3)=1$. Thus \[f_3(x_1,x_2,x_3)=\max(\min(x_1,x_2),\min(x_1,x_3),\min(x_2,x_3),\min(x_1,x_2,x_3)).\]

\bd \label{topologize1}
	Let $f=(f_1,\dots,f_n):[n]\to [n]$ be a cotransverse map. Let \[\TT(f):[0,1]^n\to [0,1]^n\] be the set map defined by 
	\[
	\TT(f)(x_1,\dots,x_n) = (\TT(f)_1(x_1,\dots,x_n),\dots,\TT(f)_n(x_1,\dots,x_n))
	\]
	with \[\TT(f)_i(x_1,\dots,x_n) = \max_{(\epsilon_1,\dots,\epsilon_n)\in f_i^{-1}(1)} \min \{x_k\mid \epsilon_k=1\}\] for all $1\leq i\leq n$.
\ed

\bp \label{restrictionf}
Let $n\geq 1$. For all $x\in [n]\subset [0,1]^n$, one has $\TT(f)(x)=f(x)$.
\ep

\bpf
It is a consequence of Proposition~\ref{maxminxformula}.
\epf

\bp \label{preserv-h-cont}
For all cotransverse maps $f:[n]\to [n]$, the set map \[\TT(f):[0,1]^n\longrightarrow [0,1]^n\] is continuous and strictly increasing. Moreover it satisfies the properties 
\[\forall (x_1,\dots,x_n)\in [0,1]^n, h(x_1,\dots,x_n) = h(\TT(f)(x_1,\dots,x_n)).\]
\ep

\bpf By Proposition~\ref{preserv-h-dis} and Proposition~\ref{restrictionf}, the theorem holds for $(x_1,\dots,x_n)\in [n] \subset [0,1]^n$. From the fact that each projection map $(x_1,\dots,x_n)\mapsto x_k$ from $[0,1]^n$ equipped with the product order to $[0,1]$ is continuous and nondecreasing, we deduce that $T(f)$ is continuous and nondecreasing. Consider a tuple $(x_1,\dots,x_n) \in [0,1]^n$. There exists a permutation $\sigma$ of $\{1,\dots,n\}$ such that $x_{\sigma(1)}\geq\dots\geq x_{\sigma(n)}$. Using Proposition~\ref{preserv-h-dis} again, write 
\begin{align*}
	& f(\epsilon_{\{\sigma(1)\}}) = \epsilon_{\{\sigma'(1)\}},\\
	& f(\epsilon_{\{\sigma(1),\sigma(2)\}}) = \epsilon_{\{\sigma'(1),\sigma'(2)\}},\\
	&\dots \\
	& f(\epsilon_{\{\sigma(1),\dots,\sigma(n)\}}) = \epsilon_{\{\sigma'(1),\dots,\sigma'(n)\}}.
\end{align*}
From the permutation $\sigma$ of $\{1,\dots,n\}$, we therefore obtain a new permutation $\sigma'$ of $\{1,\dots,n\}$. One has $\epsilon_{\{\sigma(1)\}} \in f^{-1}_{\sigma'(1)}(1)$. This means that $T(f)_{\sigma'(1)}(x_1,\dots,x_n)=x_{\sigma(1)}$ because $x_{\sigma(1)}\geq\dots \geq x_{\sigma(n)}$. One then has $\epsilon_{\{\sigma(1),\sigma(2)\}} \in f^{-1}_{\sigma'(1)}(1)$. This means that $T(f)_{\sigma'(2)}(x_1,\dots,x_n)=x_{\sigma(2)}$ because $x_{\sigma(1)}\geq\dots \geq x_{\sigma(n)}$. By repeating a finitely number of times the same argument, we obtain the equality $T(f)_{\sigma'(i)}(x_1,\dots,x_n)=x_{\sigma(i)}$ for all $1\leq i\leq n$. This implies that $T(f)(x_1,\dots,x_n)=(x_{\sigma\sigma'^{-1}(1)},\dots,x_{\sigma\sigma'^{-1}(n)})$. This means that \[h(T(f)(x_1,\dots,x_n)) = x_{\sigma\sigma'^{-1}(1)}+\dots+x_{\sigma\sigma'^{-1}(n)} = h(x_1,\dots,x_n),\] 
the first equality by definition of $h$ and the second equality since $\sigma\sigma'^{-1}$ is a permutation of $\{1,\dots,n\}$. Let $(x_1,\dots,x_n)\leq(y_1,\dots,y_n)\in [0,1]^n$. We already know that $\TT(f)(x_1,\dots,x_n) \leq \TT(f)(y_1,\dots,y_n)$. Assume that $\TT(f)(x_1,\dots,x_n) = \TT(f)(y_1,\dots,y_n)$. From the previous calculation, we obtain \[h(\TT(f)(y_1,\dots,y_n)) - h(\TT(f)(x_1,\dots,x_n)) = \sum_{i=1}^{n}(y_i-x_i) =0.\] We deduce that $(x_1,\dots,x_n)=(y_1,\dots,y_n)$. This means that $\TT(f):[0,1]^n\longrightarrow [0,1]^n$ is strictly increasing. 
\epf

Before proving Proposition~\ref{prefunc1} which leads to the definition of the cotransverse topological cube in Theorem~\ref{cotransverse-cube}, we need to establish two lemmas.

\begin{lem} \label{prefunc}
Let $f:[n]\to [n]$ and $g:[n]\to [n]$ be two cotransverse maps. Then there is the equality \[\TT(fg) = \TT(f) \TT(g).\]
\end{lem}

\bpf
Consider a tuple $(x_1,\dots,x_n) \in [0,1]^n$. We want to prove that \[\TT(fg)(x_1,\dots,x_n) = \TT(f)\TT(g)(x_1,\dots,x_n).\] Let $\sigma$ be a permutation of $\{1,\dots,n\}$ such that $x_{\sigma(1)}\geq\dots \geq x_{\sigma(n)}$. Using Proposition~\ref{preserv-h-dis}, write 
\begin{align*}
	& g(\epsilon_{\{\sigma(1)\}}) = \epsilon_{\{\sigma'(1)\}},\\
	& g(\epsilon_{\{\sigma(1),\sigma(2)\}}) = \epsilon_{\{\sigma'(1),\sigma'(2)\}},\\
	&\dots \\
	& g(\epsilon_{\{\sigma(1),\dots,\sigma(n)\}}) = \epsilon_{\{\sigma'(1),\dots,\sigma'(n)\}}
\end{align*}
for some permutation $\sigma'$ of $\{1,\dots,n\}$. From the calculation made in the proof of Proposition~\ref{preserv-h-cont}, we obtain the equality  \[(y_1,\dots,y_n)=\TT(g)(x_1,\dots,x_n)=(x_{\sigma\sigma'^{-1}(1)},\dots,x_{\sigma\sigma'^{-1}(n)}).\]
One has $y_{\sigma'(1)}\geq \dots\geq y_{\sigma'(n)}$ because $y_{\sigma'(i)} = x_{\sigma(i)}$ for all $1\leq i\leq n$. Using Proposition~\ref{preserv-h-dis} again, write 
\begin{align*}
	& f(\epsilon_{\{\sigma'(1)\}}) = \epsilon_{\{\sigma''(1)\}},\\
	& f(\epsilon_{\{\sigma'(1),\sigma'(2)\}}) = \epsilon_{\{\sigma''(1),\sigma''(2)\}},\\
	&\dots \\
	& f(\epsilon_{\{\sigma'(1),\dots,\sigma'(n)\}}) = \epsilon_{\{\sigma''(1),\dots,\sigma''(n)\}}
\end{align*}
for some permutation $\sigma''$ of $\{1,\dots,n\}$. We obtain the equality  \[\TT(f)(y_1,\dots,y_n)=(y_{\sigma'\sigma''^{-1}(1)},\dots,y_{\sigma'\sigma''^{-1}(n)}) = (x_{\sigma\sigma''^{-1}(1)},\dots,x_{\sigma\sigma''^{-1}(n)}),\]
the left-hand equality by the calculation made in the proof of Proposition~\ref{preserv-h-cont}, the right-hand equality by definition of $y_i$. Since we have 
\begin{align*}
	& fg(\epsilon_{\{\sigma(1)\}}) = f(\epsilon_{\{\sigma'(1)\}}) = \epsilon_{\{\sigma''(1)\}},\\
	& fg(\epsilon_{\{\sigma(1),\sigma(2)\}}) = f(\epsilon_{\{\sigma'(1),\sigma'(2)\}})= \epsilon_{\{\sigma''(1),\sigma''(2)\}},\\
	&\dots \\
	& fg(\epsilon_{\{\sigma(1),\dots,\sigma(n)\}}) = f(\epsilon_{\{\sigma'(1),\dots,\sigma'(n)\}})= \epsilon_{\{\sigma''(1),\dots,\sigma''(n)\}},
\end{align*}
we obtain using the calculation made in the proof of Proposition~\ref{preserv-h-cont} that \[\TT(f)\TT(g)(x_1,\dots,x_n)= \TT(f)(y_1,\dots,y_n) = \TT(fg)(x_1,\dots,x_n).\]
\epf

\begin{nota}
	For $\delta_i^\alpha:[n-1]\to [n] \in \square$, let \[\TT(\delta_i^\alpha)=
	\begin{cases}
		[0,1]^{n-1} \to [0,1]^n\\
		(\epsilon_1, \dots, \epsilon_{n-1})\mapsto (\epsilon_1,\dots, \epsilon_{i-1}, \alpha, \epsilon_i, \dots, \epsilon_{n-1})
	\end{cases}
	\] for all $n\geq 1$ and $\alpha\in\{0,1\}$.
\end{nota}

\begin{lem} \label{prefunc0.5}
Let $f:[n]\to [p]$ and $g:[m]\to [n]$ be two cotransverse maps with $f\in \square$ or $g\in \square$. Then there is the equality \[\TT(fg) = \TT(f) \TT(g).\]
\end{lem}

\bpf It is well known if both $f$ and $g$ belong to $\square$. If only one of the two maps $f$ or $g$ belongs to $\square$, we use Definition~\ref{topologize1} of $\TT(f)$ or $\TT(g)$ for the map not belonging to $\square$ and we add $0$ or $1$ to the other coordinates, depending on the coface map. 
\epf

\bp \label{prefunc1}
Let $f:[n]\to [p]$ and $g:[m]\to [n]$ be two cotransverse maps. Then there is the equality \[\TT(fg) = \TT(f) \TT(g).\]
\ep

\bpf
Consider the commutative diagram of $\widehat{\square}$ (the vertical maps are coface maps)
\[
\xymatrix@C=5em@R=3em
{
	[m] \fr{g} & [n] \fr{f} & [p]\\
	[m] \ar@{=}[u] \fr{g'} & [m]\fu{\delta} \fr{f'}& [m]\fu{\delta'}
}
\]
where the factorizations $g=\delta g'$ and $f=\delta'f'$ are given by Proposition~\ref{decomposition_distance}. Then there is the sequence of equalities (by repeatedly using Lemma~\ref{prefunc} and Lemma~\ref{prefunc0.5})
\begin{multline*}
	\TT(fg) =\TT(\delta'f'g')= \TT(\delta') \TT(f'g')= \TT(\delta') \TT(f')\TT(g')\\=\TT(\delta'f')\TT(g')=\TT(f\delta)\TT(g')
	=\TT(f)\TT(\delta)\TT(g')=\TT(f)\TT(\delta g')=\TT(f)\TT(g).
\end{multline*}
\epf

\bth \label{cotransverse-cube}
The mappings 
\begin{align*}
	&[n]\mapsto [0,1]^n \hbox{ for all } n\geq 0\\
	&f:[n]\to [n]\in \widehat{\square} \mapsto \TT(f) \hbox{ for all } n\geq 1\\
	&\delta_i^\alpha:[n-1]\to [n] \mapsto \TT(\delta_i^\alpha) \hbox{ for all } n\geq 1
\end{align*}
give rise to a cotransverse topological space called the \textit{cotransverse topological cube} and denoted by $|\widehat{\square}[*]|_{geom}$.
\eth

\bpf
The functoriality is a consequence of Proposition~\ref{prefunc1}.
\epf

Proposition~\ref{eq-cat} and Theorem~\ref{cotransverse-cube} lead to the following definition:

\bd
Let $K$ be a symmetric transverse set. Let 
\[
|K|_{geom} = \int^{[n]\in \widehat{\square}} K_n.|\widehat{\square}[n]|_{geom}
\]
This gives rise to a colimit-preserving functor $|-|_{geom}:\widehat{\square}^{op}\set \to \top$. 
\ed

A point of $|K|_{geom}$ may admit several presentations $[c;x]=|c|_{geom}(x)$ with $c\in K$ and $x\in [0,1]^{\dim (c)}$. One has $|\widehat{\square}[n]|_{geom} \iso [0,1]^n$ for all $n\geq 0$. This implies that for all cotransverse maps $f:[m]\to [n]$, by identifying using Yoneda's lemma with the map $f:\widehat{\square}[m]\to \widehat{\square}[n]$, the continuous map $|f|_{geom}:[0,1]^m\to [0,1]^n$ is the continuous map $\TT(f):[0,1]^m\to [0,1]^n$. Since all involved functors are colimit-preserving, one obtains the natural homeomorphism \[|\widehat{\mathcal{L}}(K)|_{geom}\iso |K|_{geom}\] for all precubical sets $K$ where $|K|_{geom}$ is the geometric realization of the precubical set $K$ which is defined similarly \cite[Notation~4.1]{NaturalRealization}. By Proposition~\ref{free_square}, we deduce the natural homeomorphism $|\de\widehat{\square}[n]|_{geom} \iso |\de\square[n]|_{geom}$ for all $n\geq 0$. The topology of $|K|_{geom}$ is always Hausdorff by Proposition~\ref{comp}. In the cellular case, there is a more direct proof of this fact given in Proposition~\ref{CW-complex}.

\begin{nota}
	Let $n\geq 1$. Denote by $\mathbf{D}^n = \{b\in \mathbb{R}^n, |b| \leq 1\}$ the $n$-dimensional disk, and by $\mathbf{S}^{n-1} = \{b\in \mathbb{R}^n, |b| = 1\}$ the $(n-1)$-dimensional sphere. By convention, let $\mathbf{D}^{0}=\{0\}$ and $\mathbf{S}^{-1}=\varnothing$.
\end{nota}

\bp \label{CW-complex}
	For all cellular symmetric transverse sets $K$, the geometric realization $|K|_{geom}$ is a CW-complex. In particular, the space $|K|_{geom}$ is Hausdorff.
\ep

\bpf
There are the homeomorphisms $|\widehat{\square}[n+1]|_{geom} \iso \Di^{n+1}$ and $|\de\widehat{\square}[n+1]|_{geom} \iso \Sp^n$ for $n\geq 0$ by Proposition~\ref{free_square}. Consider the diagram of solid arrows of topological spaces 
\[
\xymatrix@C=3em@R=3em
{
	|\de\widehat{\square}[n+1]|_{geom} \fd{\subset}\ar@{-->}[r]^{f} & \Sp^n \fd{\subset}\\
	|\widehat{\square}[n+1]|_{geom} \fr{\iso} & \Di^{n+1}
}
\]
Since the inclusions $|\de\widehat{\square}[n+1]|_{geom}\subset |\widehat{\square}[n+1]|_{geom}$ and $\Sp^n\subset \Di^{n+1}$ are closed inclusions, the composite map $|\de\widehat{\square}[n+1]|_{geom}\to \Di^{n+1}$ induces a homeomorphism $f:|\de\widehat{\square}[n+1]|_{geom}\to \Sp^n$ such that the diagram above is commutative. The proof is complete with Proposition~\ref{finite_induction_cofibrant}.
\epf

\begin{nota}
	Let $|\widehat{\square}[n]|_{\vd}$ be the Lawvere metric space $([0,1]^n,\vd)$ for all $n\geq 0$. 
\end{nota}

There is a topological version of Proposition~\ref{pseudoiso-dis}: 

\bp \label{pseudoiso-cont}
Let $n\geq 1$. Let $f:[n]\to [n]$ be a cotransverse map. Then $\TT(f):[0,1]^n\to [0,1]^n$ yields a map of Lawvere metric spaces from $|\widehat{\square}[n]|_{\vd}$ to itself which is quasi-isometric. 
\ep

\bpf
Let $x,y\in [0,1]^n$. Suppose first that $x \leq y$ are comparable. Then $h(x)\leq h(y)$. By Proposition~\ref{preserv-h-cont}, there is the inequality $h(\TT(f)(x)) = h(x) \leq h(y) = h(\TT(f)(y))$. We deduce that $\vd(x,y)=h(y)-h(x)=h(\TT(f)(y))-h(\TT(f)(x))=\vd(\TT(f)(x),\TT(f)(y))$, the first equality by definition of $\vd$, the second equality by the previous remark, and the last equality by definition of $\vd$ and since $\TT(f)$ is strictly increasing. Now suppose that $x\leq y$ is false. This means that $\vd(x,y)=\infty$. This implies that $\vd(\TT(f)(x),\TT(f)(y))\leq \vd(x,y)$ whatever the value of $\vd(\TT(f)(x),\TT(f)(y))$ is. Thus, $\TT(f):[0,1]^n\to [0,1]^n$ is a map of Lawvere metric spaces. 
\epf

\begin{cor} \label{quasi-isometric}
	Let $f:[m]\to [n]$ be a cotransverse map. The induced map $\TT(f):|\widehat{\square}[m]|_{\vd} \to |\widehat{\square}[n]|_{\vd}$ is a map of Lawvere metric spaces which is also quasi-isometric.
\end{cor}

This leads to the theorem: 

\bth \label{cotransverse-Lawverecube}
The mappings 
\begin{align*}
	&[n]\mapsto [0,1]^n \hbox{ for all } n\geq 0\\
	&f:[n]\to [n]\in \widehat{\square} \mapsto \TT(f) \hbox{ for all } n\geq 1\\
	&\delta_i^\alpha:[n-1]\to [n] \mapsto \TT(\delta_i^\alpha)\hbox{ for all } n\geq 1
\end{align*}
give rise to a cotransverse Lawvere metric space called the \textit{cotransverse Lawvere cube} and denoted by $|\widehat{\square}[*]|_{\vd}$.
\eth

Proposition~\ref{eq-cat} and Theorem~\ref{cotransverse-Lawverecube} lead to the following definition:

\bd \label{rea-Lawverecube}
Let $K$ be a symmetric transverse set. Let 
\[
|K|_{\vd} = \int^{[n]\in \widehat{\square}} K_n.|\widehat{\square}[n]|_{\vd}.
\]
This gives rise to a colimit-preserving functor $|-|_{\vd}:\widehat{\square}^{op}\set \to \LMet$. 
\ed

\section{Natural \mins{d}-path of a symmetric transverse set}
\label{nat-path}

It is necessary to consider the symmetric version of the notion of Lawvere metric space to obtain a convenient notion of the underlying topological space of a Lawvere metric space.

\bd
A \textit{pseudometric space} $(X,d)$ is a set $X$ equipped with a map $d:X\p X\to [0,\infty]$ called a \textit{pseudometric} such that:
\begin{itemize}
	\item $\forall x\in X,d(x,x)=0$
	\item $\forall (x,y)\in X\p X, d(x,y)=d(y,x)$ (symmetry axiom)
	\item $\forall (x,y,z)\in X\p X\p X, d(x,y)\leq d(x,z)+d(z,y)$.
\end{itemize}
A map $f:(X,d)\to (Y,d)$ of pseudometric spaces is a set map $f:X\to Y$ which is \textit{non-expansive}, i.e. $\forall (x,y)\in X\p X, d(f(x),f(y))\leq d(x,y)$.  
\ed

\begin{nota}
	The category of pseudometric spaces is denoted by $\SLMet$.
\end{nota}

The family of balls $B(x,\epsilon)=\{y\in X\mid d(x,y)<\epsilon\})$ of a pseudometric space $(X,d)$ with $x\in X$ and $\epsilon>0$ generates a topology called the \textit{underlying topology} of $(X,d)$. This construction gives rise to a functor from pseudometric spaces to general topological spaces because maps of pseudometric spaces are non-expansive. It is not colimit-preserving by \cite[Remark~3.30]{zbMATH07226006}. The category of pseudometric spaces is bicomplete, being a reflective full subcategory of the bicomplete category of Lawvere metric spaces by \cite[Proposition~3.21]{zbMATH07226006}. Start from a Lawvere metric space $(X,d)$. The image by the reflection is the pseudometric space $(X,d^\wedge)$ defined for all $(x,y)\in X\p X$ by 
\[
d^\wedge(x,y) = \min\limits_{n\geq 0}\min_{x=x_0,x_1,\dots,x_{n+2}=y}\sum_{i=0}^{n} \bigg(d(x_{i+1},x_i) + d(x_{i+1},x_{i+2})\bigg)
\]
With $(x_0,x_1,x_2)=(x,x,y)$, we obtain $d(x,y)=d(x,x)+d(x,y)\geq d^\wedge(x,y)$ for all $(x,y)\in X\p X$. 

Since there is a mistake in the statement of \cite[Proposition~3.21]{zbMATH07226006} (the formula giving $d^\wedge$ is not correct) and no proof is given, a short explanation of the adjunction is provided in this paragraph for the ease of the reader. By replacing $x_i$ by $x_{n+2-i}$ in the formula above, we deduce that $d^\wedge(x,y)=d^\wedge(y,x)$. Let us start from a map of Lawvere metric spaces $f:(X,d)\to (Y,d)$ where $(Y,d)$ is a pseudometric space. The map $f^\wedge:(X,d^\wedge)\to (Y,d)$ has the same underlying set map (so it is unique if it exists) and we just have to verify that it is non-expansive. Since $f$ is non-expansive, one has $d(f(x),f(y))\leq d(x,y)$ for all $x,y\in X$. We obtain 
\begin{multline*}
	d(f(x_0),f(x_{n+2})) \leq \bigg(\sum_{i=0}^{n-1} d(f(x_{i+1}),f(x_i))\bigg) + d(f(x_{n+1}),f(x_n)) + d(f(x_{n+1}),f(x_{n+2}))\\\leq \sum_{i=0}^{n} \bigg(d(f(x_{i+1}),f(x_i)) + d(f(x_{i+1}),f(x_{i+2}))\bigg) \leq \sum_{i=0}^{n} \bigg(d(x_{i+1},x_i) + d(x_{i+1},x_{i+2})\bigg),
\end{multline*}
the first inequality by the triangular inequality, the second inequality since one has $d(f(x_{i+1}),f(x_{i+2}))\geq 0$ for all $i\geq 0$, and the last inequality because $f:(X,d)\to (Y,d)$ is non-expansive. We deduce that $f^\wedge:(X,d^\wedge)\to (Y,d)$ is non-expansive. Conversely, if $g:(X,d^\wedge)\to (Y,d)$ is a map of pseudometric spaces, then for all $x,y\in X$, one has $d(g(x),g(y)) \leq d^\wedge(x,y)\leq d(x,y)$, the left-hand inequality since $g$ is non-expansive and the right-hand inequality by the remark above. Thus the underlying set map of $g$ induces a map of Lawvere metric spaces from $(X,d)$ to $(Y,d)$.

\bd
The \textit{underlying topological space} of a Lawvere metric space $(X,d)$ is by definition the underlying topological space of the pseudometric space $(X,d^\wedge)$. 
\ed

\begin{nota}
	Let $K$ be a symmetric transverse set. The underlying set of $|K|_{\vd}$ equipped with the pseudometric $\vd^\wedge$ gives rise to a pseudometric space denoted by $|K|_{\vd^\wedge}$. 
\end{nota}

This gives rise to a colimit-preserving functor \[|-|_{\vd^\wedge}:\widehat{\square}^{op}\set \to \LMet \to \SLMet.\] In particular, one has 
\[
|K|_{\vd^\wedge} \iso \int^{[n]\in \widehat{\square}} K_n.|\widehat{\square}[n]|_{\vd^\wedge}.
\]

\begin{nota}
	The underlying topological space of the pseudometric space $|K|_{\vd^\wedge}$ is denoted by $|K|_{d_1}$. It is a first countable (and therefore sequential) topological space, the family of balls $(B(x,1/n)_{n \geq 1})$ being a neighborhood basis of $x\in |K|_{d_1}$.
\end{nota}

\begin{nota}
	Let $n\geq 1$. Let $(x_1,\dots,x_n),(x'_1,\dots,x'_n)\in [0,1]^n$. Let
	\[
	d_1((x_1,\dots,x_n),(x'_1,\dots,x'_n)) = \sum\limits_{i=1}^{n} |x_i-x'_i|.
	\]
\end{nota}

\bp
Let $n\geq 1$. For all $x,y\in [0,1]^n$, there is the equality \[\vd^\wedge(x,y)= d_1(x,y).\] 
\ep

\bpf 
By definition, $\vd^\wedge(x,y)$ is the minimum of the sums of the form 
\begin{multline*}
	\big(\vd(x_1,x_0)+\vd(x_1,x_2)\big) + \\\big(\vd(x_2,x_1)+\vd(x_2,x_3)\big)+\dots + \big(\vd(x_{n+1},x_n)+\vd(x_{n+1},x_{n+2})\big)
\end{multline*}
with $n\geq 0$ and $x_0=x$ and $x_{n+2}=y$. To have a finite sum, the only possibility is that 
\[
x_1=x_2=\dots =x_n=x_{n+1}=z, z\leq x, z\leq y.
\]
Consequently, one has 
\[
\vd^\wedge(x,y) = \min_{\substack{z\leq x\\z\leq y}} \big(\vd(z,x) + \vd(z,y)\big) = \min_{\substack{z\leq x\\z\leq y}} \big(d_1(z,x) + d_1(z,y)\big).
\]
From the triangular inequality, we obtain $d_1(x,y)\leq \vd^\wedge(x,y)$. Write $x=x_0+x_1$ and $y=y_0+y_1$ with $x_0\leq y_0$ and $y_1\leq x_1$. Let $z=x_0+y_1$. Then one has 
\[
\vd(z,x)+\vd(z,y) = \big(h(x_1)-h(y_1)\big) + \big(h(y_0)-h(x_0)\big) = d_1(x,y).
\]
We deduce the inequality $\vd^\wedge(x,y)\leq d_1(x,y)$.
\epf

\begin{cor} \label{samecube}
	For all $n\geq 0$, there is the homeomorphism $|\widehat{\square}[n]|_{d_1} \iso [0,1]^n$.
\end{cor}

\bp \label{comp}
Let $K$ be a symmetric transverse set. Then we have the following properties:
\begin{enumerate}
	\item The underlying sets of the topological spaces $|K|_{geom}$ and $|K|_{d_1}$ are equal.
	\item The identity of the underlying set of $|K|_{geom}$ yields a continuous map from $|K|_{geom}$ to $|K|_{d_1}$.
	\item The topological spaces $|K|_{geom}$ and $|K|_{d_1}$ are Hausdorff.
	\item The topological space $|K|_{d_1}$ is $\Delta$-generated.
\end{enumerate}
\ep

\bpf Assume at first that $\top$ is the category of $\Delta$-generated spaces and let us prove the four assertions.

(1) The forgetful functor $\SLMet\to \set$ from pseudometric spaces to sets has a right adjoint given by taking a set $S$ to the pseudometric space $(S,d_0)$ with $d_0(x,y)=0$ for all $x,y\in S$. Consequently, the forgetful functor $\SLMet\to \set$ is colimit-preserving.  The forgetful functor $\top\to\set$ is topological by \cite[Proposition~3.5]{FR}, hence colimit-preserving, the category of $\Delta$-generated spaces being the final closure in the category of general topological spaces of the segment $[0,1]$. Thus, the underlying set of $|K|_{geom}$ is equal to the underlying set of $|K|_{d_1}$. 

(2) From Corollary~\ref{samecube}, we obtain the homeomorphism $|\widehat{\square}[n]|_{geom}\iso |\widehat{\square}[n]|_{d_1}$. For each $c\in K_n$, we obtain a composite continuous map
\[
\xymatrix@C=3em
{
	|\widehat{\square}[n]|_{geom}\iso |\widehat{\square}[n]|_{d_1} \fr{|c|_{d_1}} & |K|_{d_1}
}
\]
and, by the universal property of the colimit, we deduce that the identity yields a continuous map $|K|_{geom}\to |K|_{d_1}$. 

(3) It is easy to see that the pseudometric of $|K|_{\vd^\wedge}$ restricts to a metric satisfying the additional Fr\'echet axiom (i.e. the underlying topology is $T_1$) on each path-connected component thanks to the homeomorphisms $|\widehat{\square}[n]|_{d_1} \iso [0,1]^n$ for all $n\geq 0$. Thus, the topological space $|K|_{d_1}$ is Hausdorff. Since the identity map $|K|_{geom}\to |K|_{d_1}$ is one-to-one, we deduce that $|K|_{geom}$ is Hausdorff as well. 

(4) Let $x\in |K|_{d_1}$. The family of balls $(B(x,1/n)_{n \geq 1})$ is a neighborhood basis of $x$. Assume at first that $x\in K_0$. Then for all $\epsilon\in ]0,1[$, $B(x,\epsilon)$ is path-connected because each point is related to $x$ by a continuous path. Assume now that $x\in |K|_{d_1} \backslash K_0$. From the counit map $\widehat{\mathcal{L}}(\widehat{\omega}(K))\to K$ we deduce that there exists $n\geq 1$ and $c\in K_n$ such that $x=[c;(t_1,\dots,t_{\dim(c)})]$ with $(t_1,\dots,t_{\dim(c)})\in ]0,1[^{\dim(c)}$. Let $\epsilon = \min \{t_1,\dots,t_{\dim(c)}\}$. One has $\epsilon \in ]0,1[$. For all $\eta\in ]0,\epsilon[$, the ball $B(x,\eta)$ is path-connected to $x$. We have proved that the first countable topological space $|K|_{d_1}$ is locally path-connected. It is therefore $\Delta$-generated by \cite[Proposition~3.11]{MR3270173}.

Assume now that $\top$ is the category of $\Delta$-Hausdorff $\Delta$-generated spaces and let us prove the four assertions. 

(1) The inclusion functor $\top\subset \top_\Delta$ where $\top_\Delta$ is the category of $\Delta$-generated spaces has a left adjoint $w_\Delta:\top_\Delta \to \top$ by \cite[Proposition~B.7]{leftproperflow}. The topology of the colimit $|K|_{geom}$ is now given at first by taking the colimit in the category of $\Delta$-generated spaces, and then by applying $w_\Delta$ which may identify points in the underlying set. By the above proof of (2) for $\top_\Delta$, the second step is not required because the topology we obtain by taking the colimit in the category of $\Delta$-generated spaces is already Hausdorff, and therefore $\Delta$-Hausdorff. Hence the underlying sets of the topological spaces $|K|_{geom}$ and $|K|_{d_1}$ are still equal.

The proofs of (2), (3) and (4) are unchanged in this new setting.
\epf

\begin{rem}
	It is not clear whether $|K|_{geom}$ is still a CW-complex, or at least a retract of a CW-complex for any symmetric transverse set.
\end{rem}

In general, the canonical map $|K|_{geom}\to |K|_{d_1}$ of Proposition~\ref{comp} induced by the identity is not a homeomorphism, in particular for symmetric transverse sets freely generated by a locally infinite precubical set by \cite[Proposition~1.5.17]{MR1074175}. The latter proposition can be invoked because the restriction of the pseudometric $\vd^\wedge$ to each path-connected component of the topological space $|K|_{d_1}$ is a metric.

Let $U$ be a topological space. A \textit{(Moore) path} of $U$ consists of a continuous map $[0,\ell]\to U$ with $\ell>0$. Let $\gamma_1:[0,\ell_1]\to U$ and $\gamma_2:[0,\ell_2]\to U$ be two paths of a topological space $U$ such that $\gamma_1(\ell_1)=\gamma_2(0)$. The \textit{Moore composition} $\gamma_1*\gamma_2:[0,\ell_1+\ell_2]\to U$ is the Moore path defined by 
\[
(\gamma_1*\gamma_2)(t)=
\begin{cases}
	\gamma_1(t) & \hbox{ for } t\in [0,\ell_1]\\
	\gamma_2(t-\ell_1) &\hbox{ for }t\in [\ell_1,\ell_1+\ell_2].
\end{cases}
\]
The Moore composition of Moore paths is strictly associative.

\bd \label{dpath} Let $n\geq 1$. A \textit{(tame) $d$-path} of $|\widehat{\square}[n]|_{geom} = [0,1]^n$ is a nonconstant continuous map $\gamma:[0,\ell]\to [0,1]^n$ with $\ell>0$ such that $\gamma(0),\gamma(\ell)\in \{0,1\}^n$ and such that $\gamma$ is nondecreasing with respect to each axis of coordinates. Let $c\in K_n$ with $n\geq 1$ be an $n$-cube of a general symmetric transverse set $K$. A \textit{(tame) $d$-path} of $c$ is a composite continuous map denoted by $[c;\gamma]:[0,\ell] \to |K|_{geom}$ with $\ell>0$ such that $\gamma:[0,\ell]\to [0,1]^n$ is a $d$-path with $[c;\gamma]=|c|_{geom}\gamma$. Let $K$ be a general symmetric transverse set. A \textit{(tame) $d$-path} of $K$ is a continuous path $[0,\ell] \to |K|_{geom}$ which is a Moore composition $[c_1;\gamma_1] * \dots *[c_n;\gamma_n]$ of $d$-paths of the cubes $c_1,\dots,c_n$ of $K$. $\gamma(0)\in K_0$ is called the \textit{initial state} of $\gamma$ and $\gamma(\ell)\in K_0$ is called the \textit{final state} of $\gamma$. 
\ed

For all $n$-cubes $c$ of $K$ and for all cotransverse maps $f$, there are the equalities \[|f^*(c)|_{geom}=|cf|_{geom} = |c|_{geom}|f|_{geom} = |c|_{geom}\TT(f)\] by functoriality of $|-|_{geom}:\widehat{\square}^{op}\set\to \top$. This implies that there is the sequence of equalities \[[f^*(c);\gamma] = |f^*(c)|_{geom}\gamma = [c;\TT(f)\gamma]\] on $[0,\ell]$. By definition of the coend, there are also the equalities \[[c;\TT(f)\gamma(t)]=[cf;\gamma(t)] = [f^*(c);\gamma(t)]\] for all $t\in [0,\ell]$. Therefore Definition~\ref{dpath} makes sense by definition of the coend and because the continuous map $\TT(f)$ is nondecreasing by Proposition~\ref{preserv-h-cont}. 

\begin{rem}
	By convention, all $d$-paths of a symmetric transverse set $K$ start and end at a vertex of $K$.
\end{rem}

\bd \label{natural_dpath_cube} A $d$-path $\gamma=(\gamma_1,\dots,\gamma_n):[0,n]\to [0,1]^n$ of the topological $n$-cube $[0,1]^n$ is \textit{natural} if $\id_{[0,n]} = h\gamma$ (see Notation~\ref{def_h}), or more explicitly if for all $t\in [0,n]$, one has $t=\gamma_1(t)+\dots +\gamma_n(t)$. The set of natural $d$-paths of $[0,1]^n$ is denoted by $N_n$. It is equipped with the compact-open topology. 
\ed

\bp (\cite[Proposition~4.10]{NaturalRealization})
The topological space $N_n$ is $\Delta$-generated and $\Delta$-Hausdorff for all $n\geq 0$.
\ep

Another way to formulate Definition~\ref{natural_dpath_cube} is as follows: 

\bp \label{natural-metric} Equip $([0,n],\leq)$ with the Lawvere metric $\vd:[0,n]\p [0,n]\to [0,\infty]$ defined by 
\[
\vd(x,y) = \begin{cases}
	y-x & \hbox{ if }x\leq y\\
	\infty & \hbox{ if }x> y.
\end{cases}
\]
The latter Lawvere metric space is denoted by $\vec{[0,n]}$ in \cite[Example~3.2]{zbMATH07226006}. A set map $\gamma:[0,n]\to [0,1]^n$ is a natural d-path if and only if it is a quasi-isometry for $\vd$. 
\ep

\bpf
The equality $t=\gamma_1(t)+\dots +\gamma_n(t)$ for all $t\in [0,n]$ implies that any natural $d$-path is a quasi-isometry for $\vd$. Conversely, suppose that the set map $\gamma:[0,n]\to [0,1]^n$ is a quasi-isometry for $\vd$. Then by Corollary~\ref{samecube}, it is continuous for $[0,1]^n$ equipped with the standard topology. And being a quasi-isometry, it satisfies $t=\vd(0,t) = \vd(\gamma(0),\gamma(t))=\gamma_1(t)+\dots +\gamma_n(t)$ for all $t\in [0,n]$. Consequently, the continuous map $\gamma:[0,n]\to [0,1]^n$ is a natural $d$-path.
\epf

Using Proposition~\ref{natural-metric}, it is now possible to generalize to symmetric transverse sets the notion of natural $d$-path introduced by Raussen in \cite[Definition~2.14]{MR2521708} for precubical sets as follows. 

\bd \label{dpath-natural}
Let $K$ be a general symmetric transverse set. A \textit{(tame) $d$-path} of $K$ of the form $[c_1;\gamma_1] * \dots *[c_n;\gamma_n]$ is \textit{natural} if each $\gamma_i$ is a natural $d$-path of $[0,1]^{\dim(c_i)}$ in the sense of Definition~\ref{natural_dpath_cube} for $1\leq i \leq n$.
\ed 

Definition~\ref{dpath-natural} makes sense because the identity induces a continuous map from $|K|_{geom}$ to the underlying topological space $|K|_{d_1}$ of the Lawvere metric space $|K|_{\vd}$ by Proposition~\ref{comp} and because for all cotransverse maps $f$, the map $\TT(f)$ is a quasi-isometry by Corollary~\ref{quasi-isometric}.

\part{Realization of symmetric transverse sets}
\label{homotopical}

\section{The c-Reedy model structure of cotransverse objects}
\label{c-Reedy}

\bd \cite[Definition~6.12]{c-Reedy} 
Let $\C$ be a category equipped with an ordinal degree function on its objects. 
\begin{itemize}
	\item A morphism is \textit{level} if its domain and codomain have the same degree.
	\item The \textit{degree} of a factorization $(h,g)$ of a morphism $f$ is the degree of the intermediate object (i.e. the domain of $h$ which is the codomain of $g$).
	\item A factorization of a morphism $f$ is \textit{fundamental} if its degree is strictly less than the degrees of both the domain and codomain of $f$.
	\item A morphism is \textit{basic} if it does not admit any fundamental factorization. 
\end{itemize}
\ed

\bd \cite[page~37]{c-Reedy}
	Let $\C$ be a category equipped with an ordinal degree function on its objects. The \textit{$\delta$-th stratum} of $\C$, denoted by $\C_{=\delta}$, is the subcategory of $\C$ generated by the objects of degree $\delta$ and by the basic morphisms between them.
\ed

\bd
Let $\C$ be a category. Let $f$ be a map of $\C$. The \textit{category of factorizations} of $f$ has for objects the pairs $(h,g)$ such that $hg=f$ and for morphisms $k:(h,g)\to (h',g')$ the morphisms $k$ of $\C$ (which are called \textit{connecting morphisms}) such that there is a commutative diagram 
\[
\xymatrix@C=5em@R=3em
{
	\bullet \fr{g'} \ar@{=}[d] & \bullet \fr{h'}& \bullet \ar@{=}[d] \\
	\bullet \fr{g} & \fu{k}\bullet \fr{h} & \bullet 
}
\]
\ed

\bd \cite[Definition~8.25]{c-Reedy} 
A \textit{c-Reedy category} $\C$ is a small category equipped with an ordinal degree function $d$ on its objects, and subcategories $\overleftrightarrow{\C}$, $\overrightarrow{\C}$ and $\overleftarrow{\C}$ containing all objects such that 
\begin{enumerate}
	\item $\overleftrightarrow{\C}\subseteq \overrightarrow{\C} \cap\overleftarrow{\C}$.
	\item Every morphism in $\overleftrightarrow{\C}$ is level.
	\item Every morphism in $\overrightarrow{\C}\backslash \overleftrightarrow{\C}$ strictly raises degree, and every morphism in $\overleftarrow{\C}\backslash \overleftrightarrow{\C}$ strictly lowers degree.
	\item Every morphism $f$ factors as $\overrightarrow{f}\overleftarrow{f}$, where $\overrightarrow{f}\in \overrightarrow{\C}$ and $\overleftarrow{f}\in \overleftarrow{\C}$. The subcategory of the category of factorizations of $f$ generated by the pairs $(h,g)$ with $h\in\overrightarrow{\C}$ and $g\in \overleftarrow{\C}$ and such that the connecting morphisms belong to $\overleftrightarrow{\C}$ is connected for all $f$. 
	\item For any object $x$ and any degree $\delta<d(x)$, the functor $\overleftarrow{\C}(x,-):\overleftrightarrow{\C}_{=\delta} \to \set$ is a coproduct of retracts of representables.
\end{enumerate}
\ed

\begin{nota}
Let us equip the small category $\widehat{\square}$ with the ordinal degree function $d([n])=n$ for all $n\geq 0$. Let 
	\begin{align*}
		& \vec{\widehat{\square}}=\widehat{\square}\\
		& \overleftrightarrow{\widehat{\square}}=\overleftarrow{\widehat{\square}}= \coprod_{n\geq 0}\{f:[n]\to [n]\mid f\in \widehat{\square}\}
	\end{align*}
\end{nota}

\bp \label{equal}
Let $n\geq 0$. The $n$-th stratum $\widehat{\square}_{=n}$ is the full subcategory of $\widehat{\square}$ having one object $[n]$. In particular, one has \[\widehat{\square}_{=n}([n],[n])=\widehat{\square}([n],[n]).\]
\ep

\bpf
Every morphism $f:[m]\to [n]$ of $\widehat{\square}$ is basic since every factorization of $f$ as a composite $[m]\to [p]\to [n]$ implies that $m\leq p\leq n$, and therefore a factorization cannot be fundamental. Hence the proof is complete.
\epf

Let $\mathcal{M}$ be a model category. Let $\C$ be a small category. Recall that the \textit{projective model structure} is the unique model structure (if it exists) on the functor category $\mathcal{M}^\C$ such that the weak equivalences and the fibrations are the objectwise ones. 

The projective model structure on $\mathcal{M}^\C$ exists for any small category $\C$ when $\mathcal{M}$ is an accessible model category in the sense of \cite[Definition~5.1]{zbMATH06722019} or \cite[Definition~3.1.6]{HKRS17} by \cite[Theorem~3.4.1]{HKRS17} or a cofibrantly generated model category by \cite[Theorem~11.6.1]{ref_model2}.

\begin{nota}
	In this section, $\mathcal{M}$ denotes a model category such that the projective model structure on $\mathcal{M}^{\widehat{\square}_{=n}}$ exists for all $n\geq 0$. where $\widehat{\square}_{=n}$ is the $n$-th stratum.
\end{nota}

Theorem~\ref{proj-cof-suff-cond} provides a necessary and sufficient condition for a cotransverse object of $\mathcal{M}$ to be projective cofibrant; this is the analog of \cite[Proposition~2.3.1]{realization} for symmetric transverse sets. 

The key fact used in \cite{realization} is that the small category $\square$ is a direct Reedy category. This implies that the projective model structure on cocubical objects exists and that it coincides with the Reedy model structure for all model categories. It turns out that the small category $\widehat{\square}$ is not Reedy, whether understood in the sense of Berger-Moerdijk \cite[Definition~1.1]{g-Reedy} or Cisinski \cite[Definition~8.1.1]{Cisinski-Book}. 

Indeed, the factorization of a map of $\widehat{\square}$ by a map of $\overleftarrow{\widehat{\square}}$ followed by a map of $\vec{\widehat{\square}}$ is not unique up to isomorphism. For example, the following commutative diagram of $\widehat{\square}$ with $m<n$ gives rise to two non-isomorphic factorizations of $hkg$
\[
\xymatrix@C=5em@R=3em
{
	[m] \fr{kg} \ar@{=}[d] & [m] \fr{h}& [n] \ar@{=}[d] \\
	[m] \fr{g} & \fu{k}[m] \fr{hk} & [n] 
}
\]
when $k$ is non-invertible. 

However the category of such factorizations of a map has a final object by Proposition~\ref{final}. In fact, the small category $\widehat{\square}$  turns out to be c-Reedy by Proposition~\ref{example-c-Reedy}.

\bp \label{final}
Let $f$ be a map of $\widehat{\square}$. Consider the subcategory of the category of factorizations of $f$ generated by the pairs $(h,g)$ with $h\in\overrightarrow{\widehat{\square}}$ and $g\in \overleftarrow{\widehat{\square}}$. Note that the connecting morphisms are necessarily level since $\overleftrightarrow{\widehat{\square}}=\overleftarrow{\widehat{\square}}$. Then this subcategory has a final object.
\ep

\bpf
Let $f$ be a map of $\widehat{\square}$. Consider the factorization $(h,g)$ given by Proposition~\ref{decomposition_distance}: in particular, $h\in \square$. Consider another factorization $(h',g')$ of $f$. Consider the commutative diagram of solid arrows of $\widehat{\square}$
\[
\xymatrix@C=5em@R=3em
{
	\bullet \fr{g} \ar@{=}[d] & \bullet \fr{h}& \bullet \ar@{=}[d] \\
	\bullet \fr{g'} & \ar@{-->}[u]^-{k}\bullet \fr{h'} & \bullet 
}
\]
Proposition~\ref{decomposition_distance} yields the factorization $h'=h''k$ with $h''\in \square$ and $k$ which are unique. We obtain $hg=h'g'=h''kg'$. By uniqueness of the factorization of $f$ given by Proposition~\ref{decomposition_distance}, we obtain $h=h''$ and $g=kg'$, and therefore $h'=hk$. If  
\[
\xymatrix@C=5em@R=3em
{
	\bullet \fr{g} \ar@{=}[d] & \bullet \fr{h}& \bullet \ar@{=}[d] \\
	\bullet \fr{g'} & \ar@{-->}[u]^-{\overline{k}}\bullet \fr{h'} & \bullet 
}
\] 
is another commutative diagram, then $h'=hk=h\overline{k}$. By the uniqueness of Proposition~\ref{decomposition_distance}, we deduce that $k=\overline{k}$. 
\epf

\bp \label{example-c-Reedy}
The small category $\widehat{\square}$ is c-Reedy.
\ep

\bpf
One has $\overleftrightarrow{\widehat{\square}}\subset \vec{\widehat{\square}} \cap \overleftarrow{\widehat{\square}}$ (first axiom). Every morphism of $\overleftrightarrow{\widehat{\square}}$ is level (second axiom). Every morphism of $\vec{\widehat{\square}}\backslash \overleftrightarrow{\widehat{\square}}$ strictly raises degree and every morphism of $\overleftarrow{\widehat{\square}}\backslash \overleftrightarrow{\widehat{\square}}=\varnothing$ strictly lowers degree (third axiom). The category of factorizations of $f$ with connecting maps in $\overleftrightarrow{\widehat{\square}}$ is connected by Proposition~\ref{final} (fourth axiom). For every $n\geq 0$, and any degree $m<n$, the functor $\overleftarrow{\widehat{\square}}([n],-):\widehat{\square}_{=m}\to \set$ is an (empty) coproduct of retracts of representables because $\widehat{\square}([n],[m])=\varnothing$ (fifth axiom). 
\epf

\begin{nota}
	Let $\widehat{\square}_{<n}$ be the full category of $\widehat{\square}$ containing the objects $[0],\dots,[n-1]$~\footnote{This category should be denoted by $\widehat{\square}_n$ with the notation of \cite{c-Reedy}; I find this notation a bit confusing.}. 
\end{nota}

\begin{nota} \label{def}
	Let $n\geq 0$. Following the notations of \cite[page~37]{c-Reedy}, let 
	\[
	\de_n\widehat{\square}([p],[q]) = \int^{[m]\in \widehat{\square}_{<n}} \widehat{\square}([m],[q]) \p \widehat{\square}([p],[m])
	\]
The latching and matching object functors $L_n,M_n:\mathcal{M}^{\widehat{\square}}\to \mathcal{M}^{\widehat{\square}_{=n}}$ are given by
\begin{align*}
	& (M_nA)_{[n]} = \int_{[m]\in \widehat{\square}} A([m])^{\de_n\widehat{\square}([n],[m])} \\
	& (L_nA)_{[n]} = \int^{[p]\in \widehat{\square}}\de_n\widehat{\square}([p],[n]) . A([p]) 
\end{align*}
\end{nota}

We obtain: 

\bth \label{c-Reedy-model}
There exists a unique model structure on $\mathcal{M}^{\widehat{\square}}$ such that 
\begin{itemize}
	\item The weak equivalences are objectwise.
	\item A map $A\to B$ of $\mathcal{M}^{\widehat{\square}}$ is a fibration (trivial fibration resp.) if for all $n\geq 0$, the map $A([n])\to (M_nA)_{[n]}\p_{(M_nB)_{[n]}} B([n])$ is a fibration (trivial fibration resp.) of $\mathcal{M}$.
	\item A map $A\to B$ of $\mathcal{M}^{\widehat{\square}}$ is a cofibration (trivial cofibration resp.)if for all $n\geq 0$, $L_nB\sqcup_{L_nA} A\to B$ is a projective cofibration (trivial cofibration resp.) of the projective model structure of $\mathcal{M}^{\widehat{\square}_{=n}}$.
\end{itemize}
This model structure is called the \textit{c-Reedy model structure} of $\mathcal{M}^{\widehat{\square}}$. 
\eth

\bpf
By Proposition~\ref{example-c-Reedy} and \cite[Theorem~8.26]{c-Reedy}, the small category $\widehat{\square}$ is almost c-Reedy in the sense of \cite[Definition~8.8]{c-Reedy}. The proof is complete thanks to \cite[Theorem~8.9]{c-Reedy}.
\epf

\bp \label{de_calcul_func}
One has 
\[
\de_n\widehat{\square}([p],[q]) = \begin{cases}
	\varnothing & \hbox{ if } p>q \hbox{ or }n\leq p\\
	\widehat{\square}([p],[q]) & \hbox{ if } p\leq q \hbox{ and }p<n
\end{cases}
\]
\ep

\bpf
The composition induces a set map $\widehat{\square}([m],[q])\p \widehat{\square}([p],[m])  \to \widehat{\square}([p],[q])$. If $p>q$, then $\widehat{\square}([p],[q])=\varnothing$, which implies that $\widehat{\square}([m],[q]) \p \widehat{\square}([p],[m])  = \varnothing$ for all $[m]\in \widehat{\square}_{<n}$. If $n\leq p$, then $n-1<p$. This means that for all $[m]\in \widehat{\square}_n$, one has $\widehat{\square}([p],[m])=\varnothing$, which implies that $\widehat{\square}([m],[q]) \p \widehat{\square}([p],[m]) = \varnothing$ for all $[m]\in \widehat{\square}_{<n}$ as well. Assume now that $p\leq q$ and $p<n$. The set $\de_n\widehat{\square}([p],[q])$ is the quotient of 
\[
\displaystyle\coprod_{m<n} \widehat{\square}([m],[q]) \p \widehat{\square}([p],[m])
\]
by the equivalence relation generated by identifying two pairs $(h,g)$ and $(h',g')$ such that $hg=h'g'$ related by a connecting map, i.e. such that there exists a commutative diagram of $\widehat{\square}$ of the form 
\[
\xymatrix@C=5em@R=3em
{
	[p] \fr{g} \ar@{=}[d] & \bullet \fr{h}& [q] \ar@{=}[d] \\
	[p] \fr{g'} & \ar@{-->}[u]^-{k}\bullet \fr{h'} & [q] 
}
\]
Consider such a pair $(h,g)$. By applying Proposition~\ref{decomposition_distance} to $g:[p]\to [m]$, we obtain a commutative diagram of $\widehat{\square}$ of the form
\[
\xymatrix@C=5em@R=3em
{
	[p] \fr{\in \widehat{\square}} \ar@{=}[d] & [p]\fd{k\in \square} \fr{hk}& [q] \ar@{=}[d] \\
	[p] \fr{g} & [m] \fr{h} & [q] 
}
\]
This means that in $\de_n\widehat{\square}([p],[q])$, every element of $\widehat{\square}([m],[q]) \p \widehat{\square}([p],[m])$ is equivalent to an element of $\widehat{\square}([p],[q])\p \widehat{\square}([p],[p])$. Consider $(h,g)\in \widehat{\square}([p],[q]) \p \widehat{\square}([p],[p])$. By applying Proposition~\ref{decomposition_distance} to $h$, we obtain a commutative diagram of $\widehat{\square}$ of the form
\[
\xymatrix@C=5em@R=3em
{
	[p] \fr{kg} \ar@{=}[d] & [p] \fr{\in \square}& [q] \ar@{=}[d] \\
	[p] \fr{g} & \ar@{->}[u]^-{k}[p] \fr{h} & [q]
}
\]
This means that every element of $\widehat{\square}([p],[q])\p \widehat{\square}([p],[p])$ is equivalent in $\de_n\widehat{\square}([p],[q])$ to an element of $\square([p],[q])\p \widehat{\square}([p],[p])$. This means that $\de_n\widehat{\square}([p],[q])$ is the quotient of $\square([p],[q])\p \widehat{\square}([p],[p])$ by the equivalence relation. If $(h,g)$ and $(h',g')$ are two equivalent elements of $\square([p],[q]) \p \widehat{\square}([p],[p])$ in $\de_n\widehat{\square}([p],[q])$, then this implies in particular that $hg=h'g'$. By the uniqueness of the factorization given by Proposition~\ref{decomposition_distance}, this implies that $h=h'$ and $g=g'$. We obtain 
\[
\de_n\widehat{\square}([p],[q]) \iso \square([p],[q]) \p \widehat{\square}([p],[p]) \iso \widehat{\square}([p],[q]),
\]
the first isomorphism since the equivalence relation on $\square([p],[q]) \p \widehat{\square}([p],[p])$ restricts to the equality by the previous arguments, the second isomorphism by the uniqueness of Proposition~\ref{decomposition_distance}.
\epf

Let $\C$ be a small category. Consider a small diagram $X:\C\to \mathcal{M}$ and a weight $W:\C\to \set$. The weighted limit $\int_{c\in \C}X(c)^{W(c)}$ is a end which is characterized by the adjunction
\[
\mathcal{M}^{\C}(W.Y,X) \iso \mathcal{M}\bigg(Y,\int_{c\in \C}X(c)^{W(c)}\bigg).
\]
We obtain the following lemma.

\begin{lem} \label{wlim}
	Let $\C$ be a small category. Consider a small diagram $X:\C\to \mathcal{M}$ and the empty weight $W:\C\to \set$ with $W(c)=\varnothing$ for all $c\in \C$. Then there is the isomorphism 
	\[
	\int_{c\in \C}X(c)^{\varnothing} \iso \mathbf{1}.
	\]
\end{lem}

\bpf
There are the isomorphisms
\[
\mathcal{M}(Y,\mathbf{1}) \iso \mathbf{1} \iso \mathcal{M}^{\C}(\varnothing,X) \iso \mathcal{M}^{\C}(\varnothing.Y,X) \iso \mathcal{M}\bigg(Y,\int_{c\in \C}X(c)^{\varnothing}\bigg)
\]
for all objects $Y$ of $\mathcal{M}$, the right-hand isomorphism by adjunction. The proof is complete thanks to Yoneda's lemma. 
\epf

Let $\C$ be a small category. Consider a small diagram $X:\C\to \mathcal{M}$ and a weight $U:\C^{op}\to \set$. The weighted colimit $\int^{c\in \C}U(c).X(c)$ is a coend which is characterized by the adjunction
\[
\mathcal{M}^{\C}(X,Y^U) \iso \mathcal{M}\bigg(\int^{c\in \C}U(c).X(c),Y\bigg).
\]
We obtain the following lemma.

\begin{lem} \label{wcolim}
	Let $\C$ be a small category. Consider a small diagram $X:\C\to \mathcal{M}$ and a weight $U:\C^{op}\to \set$. Let $\D$ be the full subcategory of $\C$ generated by the objects $c$ such that $U(c)\neq \varnothing$. Then there is the isomorphism  
	\[
	\int^{c\in \D} U(c).X(c) \iso \int^{c\in \C} U(c).X(c).
	\]
\end{lem}

\bpf
By definition of the weighted colimits, there are the isomorphisms 
\begin{align*}
	& \mathcal{M}\bigg(\int^{c\in \D} U(c).X(c),Y\bigg) \iso \mathcal{M}^{\D}(X,Y^U)\\
	& \mathcal{M}\bigg(\int^{c\in \C} U(c).X(c),Y\bigg) \iso \mathcal{M}^{\C}(X,Y^U)
\end{align*}
for all objects $Y$ of $\mathcal{M}$. Let $\overline{\D}$ be the full subcategory of $\C$ generated by the objects $c$ such that $U(c)= \varnothing$. Let $c\in \overline{\D}$, $d\in \D$ and $f\in\C(c,d)$. Then $f$ gives rise to a set map $U(f):U(d)\to U(c)=\varnothing$, which implies that $U(d)=\varnothing$: contradiction. This means that for all $c\in \overline{\D}$ and $d\in \D$, one has $\C(c,d) =\varnothing$. By restriction, a map of $\mathcal{M}^{\C}(X,Y^U)$ gives rise to a map of $\mathcal{M}^{\D}(X,Y^U)$. Conversely, start from a map of $\mathcal{M}^{\D}(X,Y^U)$. To obtain a map of $\mathcal{M}^{\C}(X,Y^U)$, it remains to treat the case 
\[
\xymatrix@C=3em@R=3em
{
X(c) \fd{}\fr{} &(Y^U)(c) \fd{} \\
X(d) \fr{} & (Y^U)(d)
}
\]
with $c\in \C$ and $d\in \overline{\D}$. In the latter case, $(Y^U)(d)=\mathbf{1}$, which implies the natural bijection $\mathcal{M}^{\D}(X,Y^U)\iso \mathcal{M}^{\C}(X,Y^U)$. The proof is complete thanks to Yoneda's lemma. 
\epf

\bp \label{w-co-boundary}
For all $n\geq 0$, there is the isomorphism of symmetric transverse sets
\[
\de\widehat{\square}[n] \iso \int^{[p]\in \widehat{\square}_{<n}} \widehat{\square}([p],[n]).\widehat{\square}[p]
\]
\ep

\bpf
There are the isomorphisms of symmetric transverse sets
\[
\de\widehat{\square}[n] \iso \int^{[p]\in \widehat{\square}} (\de\widehat{\square}[n])_p.\widehat{\square}[p] \iso \int^{[p]\in \widehat{\square}_{<n}} (\de\widehat{\square}[n])_p.\widehat{\square}[p] \iso \int^{[p]\in \widehat{\square}_{<n}} \widehat{\square}([p],[n]).\widehat{\square}[p],
\]
the first isomorphism by applying $K=\int^{[p]\in \widehat{\square}} K_p.\widehat{\square}[p]$ to $K=\de\widehat{\square}[n]$, the second isomorphism since $(\de\widehat{\square}[n])_p=\varnothing$ for $p\geq n$ and by Lemma~\ref{wcolim}, and the last isomorphism by definition of $\de\widehat{\square}[n]$. 
\epf

\bth \label{proj-cof-suff-cond}
The projective model structure on $\mathcal{M}^{\widehat{\square}}$ exists and coincides with the c-Reedy model structure. Let $A:\widehat{\square}\to \mathcal{M}$ be a cotransverse object of $\mathcal{M}$. It is projective cofibrant if and only if for all $n\geq 0$, the map $\widehat{A}(\de\widehat{\square}[-]) \to \widehat{A}(\widehat{\square}[-])$ is a projective cofibration of $\mathcal{M}^{\widehat{\square}_{=n}}$.
\eth

The small category $\square$ is also a c-Reedy category since it is a Reedy category. In this case, there is the isomorphism of categories $\mathcal{M}\iso \mathcal{M}^{{\square}_{=n}}$ for all $n\geq 0$ and we recover \cite[Proposition~2.3.1]{realization} of the precubical setting.

\bpf
The matching object functor $M_n:\mathcal{M}^{\widehat{\square}}\to \mathcal{M}^{\widehat{\square}_{=n}}$ for all $n\geq 0$ can be calculated as follows. There is the sequence of isomorphisms of $\mathcal{M}$
\[
(M_nA)_{[n]}\iso \int_{[m]\in \widehat{\square}} A([m])^{\de_n\widehat{\square}([n],[m])} \iso \int_{[m]\in \widehat{\square}} A([m])^{\varnothing} \iso \mathbf{1},
\] 
the first isomorphism by definition of the matching object functor (Notation~\ref{def}), the second isomorphism since $\de_n\widehat{\square}([n],[m])=\varnothing$ by Proposition~\ref{de_calcul_func}, and the third isomorphism by Lemma~\ref{wlim}. Thus, the c-Reedy model structure of Theorem~\ref{c-Reedy-model} on $\mathcal{M}^{\widehat{\square}}$ coincides with the projective model structure which therefore exists. There is the sequence of isomorphisms of $\mathcal{M}$
\[
(L_nA)_{[n]} \iso \int^{[p]\in \widehat{\square}}\de_n\widehat{\square}([p],[n]) . A([p]) \iso \int^{[p]\in \widehat{\square}_{<n}}\widehat{\square}([p],[n]) . A([p]) \iso \widehat{A}(\de\widehat{\square}[n]),
\]
the first isomorphism by definition of the latching object functor (Notation~\ref{def}), the second isomorphism by Lemma~\ref{wcolim} and since $\de_n\widehat{\square}([p],[n])=\varnothing$ for $p\geq n$ by Proposition~\ref{de_calcul_func}, and finally the third isomorphism by Proposition~\ref{w-co-boundary} and since $\widehat{A}$ is colimit-preserving. By Theorem~\ref{c-Reedy-model}, the cotransverse object $A$ is projective cofibrant if and only if for all $n\geq 0$, the map $L_nA\to A$ is a projective cofibration of the projective model structure of $\mathcal{M}^{\widehat{\square}_{=n}}$. Since $A([n]) = \widehat{A}(\widehat{\square}[n])$ by definition of $\widehat{A}$, the proof is complete. 
\epf

\section{Realizing a symmetric transverse set as a flow}
\label{rea_sec}

The category $\top$ can be equipped with its q-model structure where the weak equivalences are the weak homotopy equivalences, the fibrations, called q-fibrations, are the Serre fibrations and the cofibrations, called q-cofibrations, are the retracts of relative cell complexes. The category $\top$ can also be equipped with its h-model structure thanks to \cite[Corollary~5.23]{Barthel-Riel} where the weak equivalences are the homotopy equivalences, the fibrations, called h-fibrations, are the Hurewicz fibrations and the cofibrations, called h-cofibrations, are the strong Hurewicz cofibrations. The m-model structure of $\top$ is also used in various places of the paper. The latter is obtained by mixing the q-model structure and the h-model structure using \cite[Theorem~2.1]{mixed-cole}: the weak equivalences are the weak homotopy equivalences and the fibrations are the Hurewicz fibrations. Further details are given at the very end of \cite[Appendix~B]{leftproperflow}.

\begin{nota}
	In the whole section, $r$ stands for $q$, $m$ or $h$.
\end{nota}

\bd \cite[Definition~4.11]{model3} \label{def-flow}
A \textit{flow} is a small semicategory enriched over the closed monoidal category $(\top,\p)$. The corresponding category is denoted by $\dtop$. 
\ed

A \textit{flow} $X$ consists of a topological space $\P X$ of \textit{execution paths}, a discrete space $X^0$ of \textit{states}, two continuous maps $s$ and $t$ from $\P X$ to $X^0$ called the source and target map respectively, and a continuous and associative map $*:\{(x,y)\in \P X\p \P X; t(x)=s(y)\}\longrightarrow \P X$ such that $s(x*y)=s(x)$ and $t(x*y)=t(y)$. Let $\P_{\alpha,\beta}X = \{x\in \P X\mid s(x)=\alpha \hbox{ and } t(x)=\beta\}$: it is the space of execution paths from $\alpha$ to $\beta$, $\alpha$ is called the initial state and $\beta$ is called the final state. Note that the composition is denoted by $x*y$, not by $y\circ x$. The category $\dtop$ is locally presentable by \cite[Theorem~6.11]{Moore1}.

\begin{exa}
	For a topological space $Z$, let $\glob(Z)$ be the flow defined by 
	\[
	\glob(Z)^0=\{0,1\}, \ 
	\P \glob(Z)= \P_{0,1} \glob(Z)=Z,\ 
	s=0,\  t=1.
	\]
	This flow has no composition law.
\end{exa}

\bth \label{three} \cite[Theorem~7.4]{QHMmodel} There exists a unique model structure on $\dtop$ such that: 
\begin{itemize}
	\item A map of flows $f:X\to Y$ is a weak equivalence if and only if $f^0:X^0\to Y^0$ is a bijection and for all $(\alpha,\beta)\in X^0\p X^0$, the continuous map $\P_{\alpha,\beta}X\to \P_{f(\alpha),f(\beta)}Y$ is a weak equivalence of the r-model structure of $\top$.
	\item A map of flows $f:X\to Y$ is a fibration if and only if for all $(\alpha,\beta)\in X^0\p X^0$, the continuous map $\P_{\alpha,\beta}X\to \P_{f(\alpha),f(\beta)}Y$ is a fibration of the r-model structure of $\top$.
\end{itemize}
This model structure is accessible and all objects are fibrant. It is called the r-model structure of $\dtop$.
\eth

By \cite[Theorem~7.7]{QHMmodel}, the m-model structure is the mixing of the q-model structure and the h-model structure in the sense of \cite[Theorem~2.1]{mixed-cole}. Every q-cofibration of flows is an m-cofibration and every m-cofibration of flows is an h-cofibration by \cite[Proposition~3.6]{mixed-cole}. Every h-fibration of flows is an m-fibration and every m-fibration of flows is a q-fibration by \cite[Theorem~2.1]{mixed-cole}. All involved model categories being accessible, the projective and injective r-model structures on $\dtop^\C$ exists for all small categories $\C$ by \cite[Theorem~3.4.1]{HKRS17}.

\bp \label{r-homotopy-r-cofibrant}
Let $f:X\to Y$ be a weak equivalence of the r-model structure of flows between r-cofibrant flows. Then $f$ is a weak equivalence of the h-model structure of flows.
\ep

\bpf
For $r=h$, there is nothing to prove. If $r=q$, then the spaces $\P_{\alpha,\beta}X$ and $\P_{f(\alpha),f(\beta)}Y$ are q-cofibrant by \cite[Theorem~5.7]{leftproperflow}. Using Whitehead's theorem \cite[Theorem~7.5.10, p.~124]{ref_model2}, we deduce that the map $\P_{\alpha,\beta}X\to \P_{f(\alpha),f(\beta)}Y$ is a homotopy equivalence of spaces. It remains the case $r=m$. The spaces $P_{\alpha,\beta}X$ and $\P_{f(\alpha),f(\beta)}Y$ are m-cofibrant by \cite[Theorem~8.7]{QHMmodel}. By \cite[Corollary~3.4]{mixed-cole}, we deduce that the weak homotopy equivalence $\P_{\alpha,\beta}X\to \P_{f(\alpha),f(\beta)}Y$ is a homotopy equivalence of spaces as well.
\epf

\bd \label{def-rea-flow} A functor $F:\widehat{\square}^{op}\set \to \dtop$ is an \textit{r-realization functor} (of symmetric transverse sets) if it satisfies the following properties:
\begin{itemize}
	\item $F$ is colimit-preserving.
	\item For all $n\geq 0$, the map $F(\de\widehat{\square}[n])\to F(\widehat{\square}[n])$ is an r-cofibration of $\dtop$.
	\item There is an objectwise weak equivalence of cotransverse flows $F(\widehat{\square}[*])\to \{0<1\}^*$ in the r-model structure of $\dtop$.
\end{itemize}
\ed

Theorem~\ref{exists}, Proposition~\ref{qTOmTOh} and Theorem~\ref{same-rea} prove that r-realization functors exist.

\bth \label{exists}
Let $(-)^{cof}$ be a q-cofibrant replacement functor of $\dtop$. The functor 
\[{|K|_q = \int^{[n]\in \widehat{\square}} K_n.(\{0<1\}^n)^{cof}}
\]
is a q-realization functor $|-|_q:\widehat{\square}^{op}\set \to \dtop$. 
\eth

\bpf
By \cite[Proposition~2.2.10]{symcub}, there is the isomorphism $|K|_q \iso |\widehat{\mathcal{L}}(K)|_q$ for all precubical sets $K$ where the left-hand term is the q-realization of the precubical set $K$ with the same q-cofibrant replacement functor and which is defined by 
\[{|K|_q = \int^{[n]\in {\square}} K_n.(\{0<1\}^n)^{cof}}
\]
(\cite[Theorem~3.9]{NaturalRealization}). Using Proposition~\ref{free_square}, we deduce the isomorphims of flows $|\square[n]|_q \iso |\widehat{\square}[n]|_q$ and $|\de\square[n]|_q \iso |\de\widehat{\square}[n]|_q$ for all $n\geq 0$. The proof is complete thanks to \cite[Proposition~7.4]{ccsprecub}.
\epf

\bp \label{restriction_precubical}
Let $F:\widehat{\square}^{op}\set \to \dtop$ be an r-realization functor of symmetric transverse sets. Then the composite functor $F\widehat{\mathcal{L}}:\square^{op}\set\to \dtop$ is an r-realization functor of precubical sets in the sense of \cite[Definition~3.6]{NaturalRealization}.
\ep

\bpf
A r-realization of precubical sets is a functor $G:{\square}^{op}\set \to \dtop$ which satisfies the following properties: 1) $G$ is colimit-preserving; 2) For all $n\geq 0$, the map $G(\de{\square}[n])\to G({\square}[n])$ is an r-cofibration of $\dtop$; 3) There is an objectwise weak equivalence of cotransverse flows $G(\widehat{\square}[*])\to \{0<1\}^*$ in the r-model structure of $\dtop$. The proposition is therefore a consequence of Proposition~\ref{free_square}.
\epf

\begin{rem}
	The composite functor $|\widehat{\mathcal{L}}(-)|_q:{\square}^{op}\set \to \dtop$ is the q-realization functor of precubical sets of \cite[Theorem~3.9]{NaturalRealization} with the same q-cofibrant replacement of flows.
\end{rem}

Recall that the \textit{injective model structure} is the unique model structure (if it exists) on a functor category $\mathcal{M}^I$ such that the weak equivalences and the cofibrations are the objectwise ones.

\begin{lem} (well-known~\footnote{I am unable to find a textbook expounding this elementary result.}) \label{proj-inj}
	Let $I$ be a small category. Let $\mathcal{M}$ be a model category (not necessarily cofibrantly generated) such that both the projective model structure $(\mathcal{M}^I)_{proj}$ and the injective model structure $(\mathcal{M}^I)_{inj}$ exist. Then the identity of $\mathcal{M}$ yields a left Quillen functor $(\mathcal{M}^I)_{proj}\to (\mathcal{M}^I)_{inj}$. In particular, every projective cofibration is an injective cofibration.
\end{lem}

\bpf 
Let $A\to B$ be a projective cofibration of $\mathcal{M}^I$. Let $C\to D$ be a trivial fibration of $\mathcal{M}$. Let $i\in I$. Consider a commutative diagram of $\mathcal{M}$ of the form 
\[
\xymatrix@C=3em@R=3em
{
A(i) \fr{} \fd{} & C \fd{} \\
B(i) \fr{} \ar@{-->}[ru]^-{\ell}& D
}
\]
The lift $\ell$ exists if and only, by adjunction, the lift $\overline{\ell}$ exists in the commutative diagram of $\mathcal{M}^I$
\[
\xymatrix@C=3em@R=3em
{
	A \fr{} \fd{} & C^{I(-,i)} \fd{} \\
	B \fr{} \ar@{-->}[ru]^-{\overline{\ell}}& D^{I(-,i)}
}
\]
The point is that $C^{I(-,i)}\to D^{I(-,i)}$ is a projective trivial fibration. Thus the lift $\overline{\ell}$ exists, and so does the lift $\ell$. We have proved that $A\to B$ is an injective cofibration and the proof is complete.
\epf

\bp \label{elm}
Let $F:\widehat{\square}^{op}\set \to \dtop$ be an r-realization functor. Then for all symmetric transverse sets $K$, there is a natural bijection $K_0\iso F(K)^0$. A projective r-cofibrant replacement of the cotransverse flow $F(\widehat{\square}[*])$, let us denote it by $F^{cof}(\widehat{\square}[*])$, gives rise to an r-realization functor as well.
\ep

\bpf 
From the objectwise weak equivalence of cotransverse flows $F(\widehat{\square}[*])\to \{0<1\}^*$, we deduce the objectwise bijection of cotransverse sets $F(\widehat{\square}[*])^0\iso \{0,1\}^* \iso \widehat{\square}[*]_0$. We obtain the natural bijection $F(K)^0 \iso K_0$ for all symmetric transverse sets $K$. By Theorem~\ref{proj-cof-suff-cond}, the map $F^{cof}(\de\widehat{\square}[-])\to F^{cof}(\widehat{\square}[-])$ is a projective r-cofibration of $\dtop^{\widehat{\square}_{=n}}$ for all $n\geq 0$. By Lemma~\ref{proj-inj}, the map $F^{cof}(\de\widehat{\square}[-])\to F^{cof}(\widehat{\square}[-])$ is an injective r-cofibration of $\dtop^{\widehat{\square}_{=n}}$ for all $n\geq 0$. Since the composite map $F^{cof}(\widehat{\square}[*])\longrightarrow F(\widehat{\square}[*])\longrightarrow \{0<1\}^*$ is a weak equivalence in the projective r-model structure of $\dtop^{\widehat{\square}}$, the proof is complete.
\epf

Thanks to Proposition~\ref{elm}, the following definition makes sense.

\bd \label{homotopy-r-realization}
A \textit{cofibrant r-realization} is an r-realization functor $F:\widehat{\square}^{op}\set\to \dtop$ such that the cotransverse flow $F(\widehat{\square}[*])$ is projective r-cofibrant. For an r-realization functor $F:\widehat{\square}^{op}\set\to \dtop$, the r-realization functor associated to the cotransverse flow $F^{cof}(\widehat{\square}[*])$ is called a \textit{cofibrant replacement} of $F$. It is denoted by $F^{cof}$. 
\ed

The map of cotransverse sets $F^{cof}(\widehat{\square}[*])\to F(\widehat{\square}[*])$ gives rise to a natural transformation of r-realization functor $F^{cof}\Rightarrow F$ by Proposition~\ref{eq-cat}.

\bp \label{qTOmTOh}
Every (cofibrant resp.) q-realization functor is a (cofibrant resp.) m-realization functor. Every (cofibrant resp.) m-realization functor is a (cofibrant resp.) h-realization functor.
\ep

\bpf
Let us prove at first the statements without the adjective ``cofibrant''. Every q-realization functor is an m-realization functor because every q-cofibration of flows is an m-cofibration of flows and because the weak equivalences are the same in the two model structures. Let $F:\widehat{\square}^{op}\set \to \dtop$ be an m-realization functor. Then for all $n\geq 0$, the map of flows $F(\de\widehat{\square}[n])\to F(\widehat{\square}[n])$ is an h-cofibration. The map of flows $F(\widehat{\square}[n])\to \{0<1\}^n$ is a weak equivalence of the h-model structure of flows by Proposition~\ref{r-homotopy-r-cofibrant}. We have proved that $F$ is an h-realization functor. Since every $m$-fibration of flows is a $q$-fibration, every projective q-cofibration of $\dtop^{\widehat{\square}_{=n}}$ is a projective m-cofibration of $\dtop^{\widehat{\square}_{=n}}$. Since every $h$-fibration of flows is a $m$-fibration, every projective m-cofibration of $\dtop^{\widehat{\square}_{=n}}$ is a projective h-cofibration of $\dtop^{\widehat{\square}_{=n}}$. We obtain the statements with the adjective ``cofibrant''.
\epf

\begin{cor} \label{cube_cofibrant}
	For all cofibrant symmetric transverse sets $K$ and all r-realization functors $F:\widehat{\square}^{op}\set \to \dtop$, the flow $F(K)$ is r-cofibrant. In particular, the flows $F(\de\widehat{\square}[n])$ and $F(\widehat{\square}[n])$ are r-cofibrant for all $n\geq 0$. 
\end{cor}

\bpf
Since $K$ is cofibrant and $F$ colimit-preserving, the map $\varnothing \to F(K)$ is a retract of a transfinite composition of pushouts of maps of the form  $F(\de\widehat{\square}[n])\to F(\widehat{\square}[n])$ for $n\geq 0$. Thus, $F(K)$ is r-cofibrant. The second statement is a consequence of Proposition~\ref{cp_cofibrant}. 
\epf

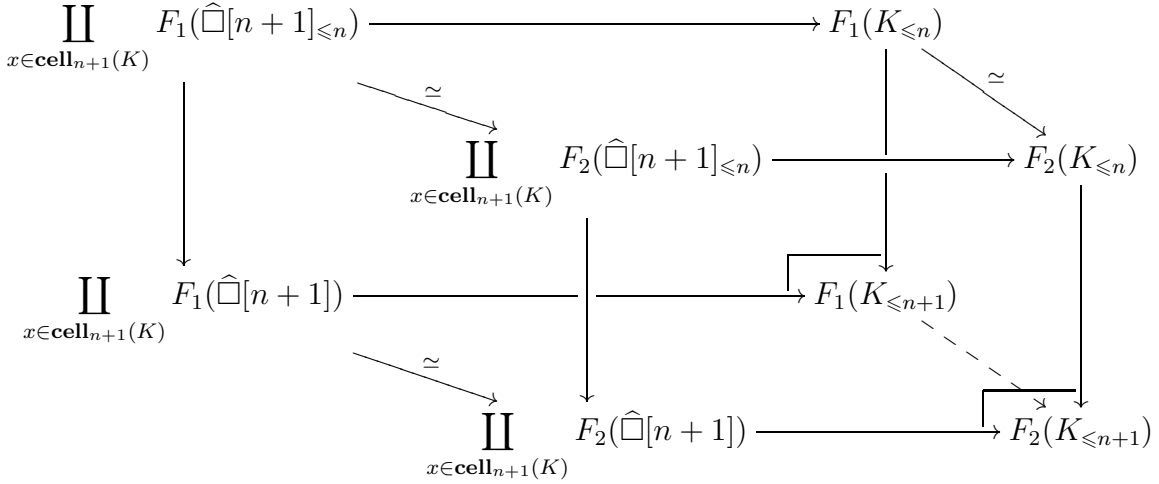
\begin{figure}
	\[{
		\xymatrix@C=1em@R=1.5em{ \displaystyle\coprod\limits_{x\in \cell_{n+1}(K)} F_1(\widehat{\square}[n+1]_{\leq n})
			\ar@{->}[rd]^-{\simeq} \ar@{->}[dd] \ar@{->}[rr] && F_1(K_{\leq
				n})
			\ar@{->}[rd]^-{\simeq} \ar@{->}'[d][dd]  &\\
			& \displaystyle\coprod\limits_{x\in \cell_{n+1}(K)}
			F_2(\widehat{\square}[n+1]_{\leq n}) \ar@{->}[rr] \ar@{->}[dd] && F_2(K_{\leq n}) \ar@{->}[dd] \\
			\displaystyle\coprod\limits_{x\in \cell_{n+1}(K)}
			F_1(\widehat{\square}[n+1]) \ar@{->}[rd]^-{\simeq} \ar@{->}'[r][rr] &&
			\cocartesien F_1(K_{\leq n+1}) \ar@{-->}[rd]^-{}& \\
			& \displaystyle\coprod\limits_{x\in \cell_{n+1}(K)} F_2(\widehat{\square}[n+1]) \ar@{->}[rr] &&
			\cocartesien F_2(K_{\leq n+1})}
	} \]
	\caption{From $n$ to $n+1$}
	\label{passage}
\end{figure}

\bp \label{elm2}
Let $F_1,F_2:\widehat{\square}^{op}\set \to \dtop$ be two r-realization functors. Suppose that there exists a commutative diagram of cotransverse flows 
\[
\xymatrix@C=3em@R=3em
{
	F_1(\widehat{\square}[*]) \fr{} \fd{} & F_2(\widehat{\square}[*]) \fd{} \\
	\{0<1\}^* \ar@{=}[r] & \{0<1\}^*
}
\]
Then the above hypothesis yields a natural map of flows $F_1(K)\to F_2(K)$ for all symmetric transverse sets $K$ which is, for all cofibrant symmetric transverse sets $K$, a weak equivalence of the r-model structure of $\dtop$ between r-cofibrant flows. Moreover, for all $(\alpha,\beta)\in K_0\p K_0$, the continuous map $\P_{\alpha,\beta}F_1(K)\to \P_{\alpha,\beta}F_2(K)$ is a homotopy equivalence between r-cofibrant topological spaces for all cofibrant symmetric transverse sets $K$. 
\ep

\bpf
By the \ttt, for all $n\geq 0$, the map of flows $F_1(\widehat{\square}[n])\to F_2(\widehat{\square}[n])$ is a weak equivalence of the r-model structure of $\dtop$, and moreover between r-cofibrant flows by Corollary~\ref{cube_cofibrant}. By Proposition~\ref{eq-cat}, the hypotheses of the proposition yield a natural transformation $\mu:F_1\Rightarrow F_2$. Let us prove by induction on $n\geq 0$ that the canonical map $F_1(K_{\leq n}) \to F_2(K_{\leq n})$ is a weak equivalence of the r-model structure of flows between r-cofibrant flows for all cellular symmetric transverse sets $K$. By Proposition~\ref{elm}, there are the natural bijections $F_i(K_{\leq 0}) = K_0\iso F_i(K)^0$ for $i=1,2$. Thus the induction hypothesis is proved for $n=0$. Let $n\geq 0$. Using the existence of the natural transformation $F_1\Rightarrow F_2$ and thanks to Proposition~\ref{finite_induction_cofibrant}, the passage from $n$ to $n+1$ can be depicted by the diagram of flows of Figure~\ref{passage}. By the induction hypothesis, and since $\widehat{\square}[n+1]_{\leq n}$ is cellular by Proposition~\ref{cp_cofibrant}, the maps of flows $F_1(\widehat{\square}[n+1]_{\leq n}) \to F_2(\widehat{\square}[n+1]_{\leq n})$ and $F_1(K_{\leq n}) \to F_2(K_{\leq n})$ are weak equivalences of the r-model structure of flows between r-cofibrant flows. We have already seen above that the map of flows $F_1(\widehat{\square}[n+1]) \to F_2(\widehat{\square}[n+1])$ is also a weak equivalence of the r-model structure of flows between r-cofibrant flows. By definition of an r-realization functor, we can apply the cube lemma \cite[Proposition~15.10.10]{ref_model2} \cite[Lemma~5.2.6]{MR99h:55031} in the r-model structure of $\dtop$ to conclude that the map $F_1(K_{\leq n+1})\to F_2(K_{\leq n+1})$ is a weak equivalence of the r-model structure of $\dtop$ between r-cofibrant flows. Since the colimits $\liminj F_1(K_{\leq n})$ and $\liminj F_2(K_{\leq n})$ are colimits of towers of r-cofibrations between r-cofibrant flows, they are homotopy colimits by \cite[Proposition~15.10.12]{ref_model2}. We conclude that the map of flows $F_1(K)\to F_2(K)$ is a weak equivalence of the r-model structure of $\dtop$ between r-cofibrant flows for all cellular symmetric transverse sets $K$. We deduce the same assertion for all cofibrant symmetric transverse sets $K$. The proof is complete thanks to Proposition~\ref{r-homotopy-r-cofibrant}. 
\epf

Proposition~\ref{elm2} has two corollaries.

\begin{cor} \label{q1}
	Let $F:\widehat{\square}^{op}\set \to \dtop$ be an r-realization functor. Then for all cofibrant symmetric transverse sets $K$, the map $F^{cof}(K)\to F(K)$ is a weak equivalence of the r-model structure of $\dtop$ between r-cofibrant flows. Moreover, for all $(\alpha,\beta)\in K_0\p K_0$, the continuous map $\P_{\alpha,\beta}F^{cof}(K)\to \P_{\alpha,\beta}F(K)$ is a homotopy equivalence between r-cofibrant topological spaces for all cofibrant symmetric transverse sets $K$.
\end{cor}

By \cite[Proposition~1.3]{MR1780498}, there exists a (non unique) functorial factorization $\varnothing \to K^{cof}\to K$ by an element of $\cell(\{\de\widehat{\square}[n]\to \widehat{\square}[n]\mid n\geq 0\})$ followed by an element of $\inj(\{\de\widehat{\square}[n]\to \widehat{\square}[n]\mid n\geq 0\})$. The functor $(-)^{cof}:\widehat{\square}^{op}\set \to \widehat{\square}^{op}\set$ is called a \textit{cofibrant replacement} of $K$. Corollary~\ref{q2} is a reformulation of Corollary~\ref{q1}.

\begin{cor} \label{q2}
	Let $F:\widehat{\square}^{op}\set \to \dtop$ be an r-realization functor. Then for all symmetric transverse sets $K$, the map $F^{cof}(K^{cof})\to F(K^{cof})$ is a weak equivalence of the r-model structure of $\dtop$ between r-cofibrant flows. Moreover, for all $(\alpha,\beta)\in K_0\p K_0$, the continuous map $\P_{\alpha,\beta}F^{cof}(K^{cof})\to \P_{\alpha,\beta}F(K^{cof})$ is a homotopy equivalence between r-cofibrant topological spaces for all symmetric transverse sets $K$.
\end{cor}

\bth \label{homotopy-natural}
Consider two r-realization functors \[F_1,F_2:\widehat{\square}^{op}\set \longrightarrow \dtop.\] Then there exists a cofibrant r-realization functor $F_3$ and a zigzag of natural transformations \[F_1\Longleftarrow F_3 \Longrightarrow F_2\] such that there is a commutative diagram of cotransverse flows 
\[
\xymatrix@C=3em@R=3em
{
	F_1(\widehat{\square}[*])\fd{} &\fl{} F_3(\widehat{\square}[*]) \fd{}\fr{} & F_2(\widehat{\square}[*]) \fd{} \\
	{\{0<1\}^*} \ar@{=}[r] & {\{0<1\}^*} \ar@{=}[r] & {\{0<1\}^*}
}
\]
and such that for all cofibrant symmetric transverse sets $K$, the maps $F_3(K)\to F_1(K)$ and $F_3(K)\to F_2(K)$ natural with respect to $K$ are weak equivalences of the r-model structure of $\dtop$ between r-cofibrant flows. Moreover, for all $(\alpha,\beta)\in K_0\p K_0$, the natural maps $\P_{\alpha,\beta}F_3(K)\stackrel{\simeq}\longrightarrow \P_{\alpha,\beta}F_1(K)$ and $\P_{\alpha,\beta}F_3(K)\stackrel{\simeq}\longrightarrow \P_{\alpha,\beta}F_2(K)$ are homotopy equivalences between r-cofibrant topological spaces for all cofibrant symmetric transverse sets $K$. When e.g. $F_1$ is already cofibrant as an r-realization functor, one can suppose that $F_1=F_3$.
\eth

\bpf
Let $F_3=F_1^{cof}$. Consider the diagram of solid arrows of $\dtop^{\widehat{\square}}$
\[
\xymatrix@C=3em@R=3em
{
	F_3(\widehat{\square}[*]) \ar@{-->}[r]^-{\mu} \fd{} & F_2(\widehat{\square}[*]) \fd{} \\
	F_1(\widehat{\square}[*]) \fr{} & \{0<1\}^*
}
\]
Since all spaces of execution paths of $\{0<1\}^*$ are discrete, the right vertical map is a trivial projective r-fibration of $\dtop^{\widehat{\square}}$. Thus, there exists a map of cotransverse flows $\mu:F_3(\widehat{\square}[*])\to F_2(\widehat{\square}[*])$ making commutative the diagram above. The proof is complete thanks to Proposition~\ref{eq-cat} and Proposition~\ref{elm2}. Assume now that $F_1$ is already cofibrant. Consider the diagram of solid arrows of $\dtop^{\widehat{\square}}$
\[
\xymatrix@C=3em@R=3em
{
	 & F_2(\widehat{\square}[*]) \fd{} \\
	F_1(\widehat{\square}[*]) \ar@{-->}[ru]^-{\mu}\fr{} & \{0<1\}^*
}
\]
The cotransverse flow $F_1(\widehat{\square}[*])$ is projective r-cofibrant. The vertical map is a trivial projective r-fibration of $\dtop^{\widehat{\square}}$. Hence the proof is complete.
\epf

Note the difference with the precubical case of \cite[Theorem~3.8]{NaturalRealization}. There is, in general, in the symmetric transverse setting, no natural transformation from $F_1$ to $F_2$. The point is that, in the precubical setting, the category $\square_{=n}$ (see Proposition~\ref{equal}) is the terminal category for all $n\geq 0$. Thus, in the precubical setting, every r-realization functor $F$ corresponds to a projective r-cofibrant cocubical flow $F(\square[*])$. In fact there is the proposition: 

\bth \label{cofibrant-case}
Consider two cofibrant r-realization functors \[F_1,F_2:\widehat{\square}^{op}\set \longrightarrow \dtop.\] Then there exists a natural transformation $F_1\Rightarrow F_2$ such that there is a commutative diagram of cotransverse flows 
\[
\xymatrix@C=3em@R=3em
{
	 F_1(\widehat{\square}[*]) \fd{}\fr{} & F_2(\widehat{\square}[*]) \fd{} \\
	 {\{0<1\}^*} \ar@{=}[r] & {\{0<1\}^*}
}
\]
and such that for all symmetric transverse sets $K$ (not necessarily cofibrant), the map $F_1(K)\to F_2(K)$ natural with respect to $K$ is a weak equivalence of the r-model structure of $\dtop$ between r-cofibrant flows. Moreover, for all $(\alpha,\beta)\in K_0\p K_0$, the natural map $\P_{\alpha,\beta}F_1(K)\stackrel{\simeq}\longrightarrow \P_{\alpha,\beta}F_2(K)$ is a homotopy equivalence between r-cofibrant topological spaces for all symmetric transverse sets $K$ (not necessarily cofibrant). 
\eth

\bpf
The existence of the natural transformation is given by Theorem~\ref{homotopy-natural}. Let $K$ be a symmetric transverse set. Consider the comma category $(\widehat{\square}\ddownarrow K)$ whose objects are the maps of transverse sets $\widehat{\square}[n]\to K$ and whose maps are the commutative squares 
\[
\xymatrix@C=3em@R=3em
{
\widehat{\square}[m]\fr{} \fd{} &  K \ar@{=}[d] \\
\widehat{\square}[n]\fr{} &  K
}
\]
We adapt Proposition~\ref{decomposition_distance} to the comma category $(\widehat{\square}\ddownarrow K)$ as follows. Since there is the equality  $\widehat{\square}^{op}\set(\widehat{\square}[p],\widehat{\square}[q]) = \widehat{\square}([p],[q])$ for all $p,q\geq 0$ by Yoneda's lemma, a map \[\widehat{\square}[m]\longrightarrow \widehat{\square}[n] \longrightarrow K\] factors uniquely as a composite map \[\widehat{\square}[m]\longrightarrow \widehat{\square}[m]\longrightarrow \widehat{\square}[n] \longrightarrow K\] such that the map $\widehat{\square}[m]\to \widehat{\square}[n]$ corresponds to a coface map by Yoneda's lemma. We obtain that  the comma category $(\widehat{\square}\ddownarrow K)$ is a c-Reedy category by mimicking the proof of Proposition~\ref{example-c-Reedy}. From the isomorphisms of symmetric transverse sets (symmetric transverse sets being presheaves over $\widehat{\square}$)
\[
K\iso \int^{[n]\in \widehat{\square}} K_n.\widehat{\square}[n] \iso \liminj_{\widehat{\square}[n]\to K} \widehat{\square}[n],
\]
we deduce for $i=1,2$ the isomorphisms of flows
\[
F_i(K) \iso \int^{[n]\in \widehat{\square}} K_n.F_i(\widehat{\square}[n]) \iso \liminj_{\widehat{\square}[n]\to K} F_i(\widehat{\square}[n])
\]
since the functor $F_i$ is colimit-preserving. Since $F_i$ is a cofibrant r-realization by hypothesis, we deduce that the right-hand colimit is a homotopy colimit in the r-model structure of flows by adapting the proof of Theorem~\ref{proj-cof-suff-cond} to $\dtop^{(\widehat{\square}\ddownarrow K)}$. By the \ttt, the natural map \[F_1(\widehat{\square}[n])\longrightarrow F_2(\widehat{\square}[n])\] is a weak equivalence of the r-model structure of $\dtop$, and moreover between r-cofibrant flows by Corollary~\ref{cube_cofibrant}. Hence the map \[F_1(K)\longrightarrow F_2(K)\] natural with respect to $K$ is a weak equivalences of the r-model structure of $\dtop$ between r-cofibrant flows for all symmetric transverse sets $K$. The proof is complete thanks to Proposition~\ref{r-homotopy-r-cofibrant}. 
\epf

\begin{qu}
	For all symmetric transverse sets $K$, there is a natural map of flows $F^{cof}(K^{cof})\to F^{cof}(K)$ for all r-realization functors $F:\widehat{\square}^{op}\set\to \dtop$. Is this natural map a weak equivalence of the r-model structure of flows ?
\end{qu}

\section{Natural realization of a symmetric transverse set}
\label{concl}

We want to use the notion of natural $d$-path of a symmetric transverse set introduced in Section~\ref{nat-path} to build the natural realization functor from symmetric transverse sets to flows, exactly as we proceed in \cite[Section~5]{NaturalRealization} for precubical sets. The definition is almost a copy-pasting. However, the verification of the functoriality is a little bit more complicated than in the precubical setting: see Proposition~\ref{incl}.

We define a flow $|\widehat{\square}[n]|_{nat}$ for $n\geq 0$ called \textit{the natural $n$-cube} as follows. The set of states is $\{0,1\}^n$. Let $n\geq 1$ and $\alpha,\beta\in\{0,1\}^n$. Recall that the topological space $N_m$ of natural $d$-paths of $[0,1]^m$ for $m\geq 1$ is defined in Definition~\ref{natural_dpath_cube}. Let
\[
\P_{\alpha,\beta}|\widehat{\square}[n]|_{nat} = \begin{cases}
	N_m & \hbox{ if }\vd(\alpha,\beta) = m\geq 1 \hbox{ and }\alpha<\beta\\
	\varnothing & \hbox{ if }\alpha\geq \beta.
\end{cases}
\]
The map $[0,1]^{m_1}\sqcup [0,1]^{m_2} \to [0,1]^{m_1+m_2}$ defined by taking $(t_1,\dots,t_{m_1})$ to $(t_1,\dots,t_{m_1},0_{m_2})$ and $(t'_1,\dots,t'_{m_2})$ to $(1_{m_1},t'_1,\dots,t'_{m_2})$ induces a continuous map $N_{m_1}\p N_{m_2}\to N_{m_1+m_2}$ by using the fact that the Moore composition of two natural $d$-paths is still a natural $d$-path. It yields the associative composition law of the flow $|\widehat{\square}[n]|_{nat}$.

Let $f:[m]\to [n]$ be a cotransverse map. Let $\alpha,\beta\in\{0,1\}^m$. Assume that $k=\vd(\alpha,\beta) \geq 1$.  There exists a unique coface map $\delta:[k]\to [m]$ with takes $0_k$ to $\alpha$ and $1_k$ to $\beta$. Consider the commutative diagram of $\widehat{\square}$ where the vertical maps are coface maps 
\[
\xymatrix@C=4em@R=4em
{
[m] \fr{f} & [n]\\
[k] \fu{\delta} \fr{[f]_{\alpha,\beta}} & [k]\fu{\delta'}
}
\]
obtained by applying Proposition~\ref{decomposition_distance} to $f\delta$. Then the continuous map $\P_{\alpha,\beta}|\widehat{\square}[m]|_{nat}\to \P_{f(\alpha),f(\beta)}|\widehat{\square}[n]|_{nat}$ induced by $f$ is the continuous map $\TT([f]_{\alpha,\beta}):N_k\to N_k$.

\bp \label{incl}
	We obtain a well-defined cotransverse flow $|\widehat{\square}[*]|_{nat}$. 
\ep

\bpf
Let $f:[m]\to [n]$ and $g:[n]\to [p]$ be two cotransverse maps. Let $\alpha,\beta\in\{0,1\}^m$. Assume that $k=\vd(\alpha,\beta) \geq 1$. Consider the commutative diagram of $\widehat{\square}$ where the vertical maps are coface maps:
\[
\xymatrix@C=5em@R=3em
{
	[m] \fr{f} & [n] \fr{g} & [p]\\
	[k] \fu{\delta} \fr{[f]_{\alpha,\beta}} & [k]\fu{\delta'} \fr{[g]_{f(\alpha),f(\beta)}}& [k]\fu{\delta''}
}
\]
Because of the uniqueness of the factorization given by Proposition~\ref{decomposition_distance}, we have
\[
[gf]_{\alpha,\beta} = [g]_{f(\alpha),f(\beta)}[f]_{\alpha,\beta}.
\]
We obtain 
\[
\TT([gf]_{\alpha,\beta}) = \TT([g]_{f(\alpha),f(\beta)})\TT([f]_{\alpha,\beta})
\]
by Proposition~\ref{prefunc1}.
\epf

Using Proposition~\ref{eq-cat} and Proposition~\ref{incl}, we obtain:

\bd \label{natural-rea}
Let $K$ be a symmetric transverse set. Consider the colimit-preserving functor 
\[
|K|_{nat} = \int^{[n]\in \widehat{\square}} K_n. |\widehat{\square}[n]|_{nat}.
\]
It is called the \textit{natural realization of $K$} as a flow.
\ed

\bp \label{extend_nat}
The composite functor $|\widehat{\mathcal{L}}(-)|_{nat}:{\square}^{op}\set \to \dtop$ (cf. Notation~\ref{free-transverse}) is the natural realization functor of precubical sets of \cite[Definition~5.3]{NaturalRealization}.
\ep

\bpf
One has $|\square[n]|_{nat} = |\widehat{\square}[n]|_{nat}$: the natural realization of the precubical set $\square[n]$ is equal to the natural realization of the symmetric transverse set $\widehat{\square}[n]$ because it is exactly the same definition. Using Proposition~\ref{free_square}, we deduce for all $n\geq 0$ the natural isomorphism $|\square[n]|_{nat} \iso |\widehat{\mathcal{L}}(\square[n])|_{nat}$. Since all involved functors are colimit-preserving, we obtain for all precubical sets $K$ the isomorphism of flows $|K|_{nat} \iso |\widehat{\mathcal{L}}(K)|_{nat}$. 
\epf

The following theorem concludes the paper.

\bth \label{same-rea}
The natural realization functor $|-|_{nat}$ from symmetric transverse sets to flows defined in Definition~\ref{natural-rea} is an m-realization functor. Let $|-|_{q}:\widehat{\square}^{op}\set \to \dtop$ be a q-realization functor. There exists an m-realization functor $F:\widehat{\square}^{op}\set \to \dtop$ and two natural transformations inducing bijections on the sets of states \[|-|_q\Longleftarrow F(-)\Longrightarrow |-|_{nat}\] such that for all cofibrant symmetric transverse sets $K$ and all $(\alpha,\beta)\in K_0\p K_0$, there is the zigzag of natural homotopy equivalences between m-cofibrant topological spaces 
\[\xymatrix
{
\P_{\alpha,\beta} |K|_q & \fl{\simeq}  \P_{\alpha,\beta} F(K) \fr{\simeq} &\P_{\alpha,\beta} |K|_{nat}
}
.\] If $|-|_{q}$ is cofibrant as a q-realization functor, then one can suppose that $F= |-|_{q}$.
\eth

\bpf
Using Proposition~\ref{extend_nat} and Proposition~\ref{free_square}, we obtain the isomorphism of flows $|\de\square[n]|_{nat} \iso |\de\widehat{\square}[n]|_{nat}$ for all $n\geq 0$. Thus the natural realization functor from symmetric transverse sets to flow is an m-realization functor because the natural realization functor of precubical sets as a flow is an m-realization functor by \cite[Theorem~5.9]{NaturalRealization}. Every q-realization functor is an m-realization functor by Proposition~\ref{qTOmTOh}. The proof is complete thanks to Theorem~\ref{homotopy-natural}.
\epf

\part*{References}

\renewcommand{\section}[2]{}%

\begin{thebibliography}{10}
	
	\bibitem{TheBook}
	J.~Ad{\'a}mek and J.~Rosick{\'y}.
	\newblock {\em Locally presentable and accessible categories}.
	\newblock Cambridge University Press, Cambridge, 1994.
	\newblock \href {https://doi.org/10.1017/cbo9780511600579.004}
	{\path{https://doi.org/10.1017/cbo9780511600579.004}}.
	
	\bibitem{Barthel-Riel}
	T.~Barthel and E.~Riehl.
	\newblock On the construction of functorial factorizations for model
	categories.
	\newblock {\em Algebr. Geom. Topol.}, 13(2):1089--1124, 2013.
	\newblock \href {https://doi.org/10.2140/agt.2013.13.1089}
	{\path{https://doi.org/10.2140/agt.2013.13.1089}}.
	
	\bibitem{MR1780498}
	T.~Beke.
	\newblock Sheafifiable homotopy model categories.
	\newblock {\em Math. Proc. Cambridge Philos. Soc.}, 129(3):447--475, 2000.
	\newblock \href {https://doi.org/10.1017/S0305004100004722}
	{\path{https://doi.org/10.1017/S0305004100004722}}.
	
	\bibitem{g-Reedy}
	C.~Berger and I.~Moerdijk.
	\newblock On an extension of the notion of reedy category.
	\newblock {\em Math. Z.}, 269(3-4):977--1004, September 2010.
	\newblock \href {https://doi.org/10.1007/s00209-010-0770-x}
	{\path{https://doi.org/10.1007/s00209-010-0770-x}}.
	
	\bibitem{Brown_cube}
	R.~Brown and P.~J. Higgins.
	\newblock On the algebra of cubes.
	\newblock {\em J. Pure Appl. Algebra}, 21(3):233--260, 1981.
	\newblock \href {https://doi.org/10.1016/0022-4049(81)90018-9}
	{\path{https://doi.org/10.1016/0022-4049(81)90018-9}}.
	
	\bibitem{MR3270173}
	D.~Christensen, G.~Sinnamon, and E.~Wu.
	\newblock The {$D$}-topology for diffeological spaces.
	\newblock {\em Pacific J. Math.}, 272(1):87--110, 2014.
	\newblock \href {https://doi.org/10.2140/pjm.2014.272.87}
	{\path{https://doi.org/10.2140/pjm.2014.272.87}}.
	
	\bibitem{Cisinski-Book}
	D-C. Cisinski.
	\newblock {\em Les pr{\'e}faisceaux comme mod{\`e}les des types d'homotopie},
	volume 308 of {\em Ast{\'e}risque}.
	\newblock Paris: Soci{\'e}t{\'e} Math{\'e}matique de France, 2006.
	
	\bibitem{mixed-cole}
	M.~Cole.
	\newblock Mixing model structures.
	\newblock {\em Topology Appl.}, 153(7):1016--1032, 2006.
	\newblock \href {https://doi.org/10.1016/j.topol.2005.02.004}
	{\path{https://doi.org/10.1016/j.topol.2005.02.004}}.
	
	\bibitem{DAT_book}
	L.~Fajstrup, E.~Goubault, E.~Haucourt, S.~Mimram, and M.~Raussen.
	\newblock {\em Directed algebraic topology and concurrency. {With} a foreword
		by {Maurice} {Herlihy} and a preface by {Samuel} {Mimram}}.
	\newblock SpringerBriefs Appl. Sci. Technol. Springer, 2016.
	\newblock \href {https://doi.org/10.1007/978-3-319-15398-8}
	{\path{https://doi.org/10.1007/978-3-319-15398-8}}.
	
	\bibitem{FR}
	L.~Fajstrup and J.~Rosick{\'y}.
	\newblock A convenient category for directed homotopy.
	\newblock {\em Theory Appl. Categ.}, 21:7--20, 2008.
	
	\bibitem{MR1074175}
	R.~Fritsch and R.~A. Piccinini.
	\newblock {\em Cellular structures in topology}, volume~19 of {\em Cambridge
		Studies in Advanced Mathematics}.
	\newblock Cambridge University Press, Cambridge, 1990.
	\newblock \href {https://doi.org/10.1017/cbo9780511983948}
	{\path{https://doi.org/10.1017/cbo9780511983948}}.
	
	\bibitem{model3}
	P.~Gaucher.
	\newblock A model category for the homotopy theory of concurrency.
	\newblock {\em Homology Homotopy Appl.}, 5(1):p.549--599, 2003.
	\newblock \href {https://doi.org/10.4310/hha.2003.v5.n1.a20}
	{\path{https://doi.org/10.4310/hha.2003.v5.n1.a20}}.
	
	\bibitem{realization}
	P.~{Gaucher}.
	\newblock {Globular realization and cubical underlying homotopy type of time
		flow of process algebra}.
	\newblock {\em {New York J. Math.}}, 14:101--137, 2008.
	
	\bibitem{ccsprecub}
	P.~Gaucher.
	\newblock Towards a homotopy theory of process algebra.
	\newblock {\em Homology Homotopy Appl.}, 10(1):353--388, 2008.
	\newblock \href {https://doi.org/10.4310/HHA.2008.v10.n1.a16}
	{\path{https://doi.org/10.4310/HHA.2008.v10.n1.a16}}.
	
	\bibitem{symcub}
	P.~Gaucher.
	\newblock Combinatorics of labelling in higher-dimensional automata.
	\newblock {\em Theor. Comput. Sci.}, 411(11-13):1452--1483, March 2010.
	\newblock \href {https://doi.org/10.1016/j.tcs.2009.11.013}
	{\path{https://doi.org/10.1016/j.tcs.2009.11.013}}.
	
	\bibitem{pastsim}
	P.~{Gaucher}.
	\newblock {Combinatorics of past-similarity in higher dimensional transition
		systems}.
	\newblock {\em {Theory Appl. Categ.}}, 32:1107--1164, 2017.
	
	\bibitem{Moore1}
	P.~Gaucher.
	\newblock Homotopy theory of {M}oore flows ({I}).
	\newblock {\em {Compositionality}}, 3(3), 2021.
	\newblock \href {https://doi.org/10.32408/compositionality-3-3}
	{\path{https://doi.org/10.32408/compositionality-3-3}}.
	
	\bibitem{leftproperflow}
	P.~Gaucher.
	\newblock Left properness of flows.
	\newblock {\em Theory Appl. Categ.}, 37(19):562--612, 2021.
	
	\bibitem{QHMmodel}
	P.~Gaucher.
	\newblock Six model categories for directed homotopy.
	\newblock {\em Categ. Gen. Algebr. Struct. Appl.}, 15(1):145--181, 2021.
	\newblock \href {https://doi.org/10.52547/cgasa.15.1.145}
	{\path{https://doi.org/10.52547/cgasa.15.1.145}}.
	
	\bibitem{NaturalRealization}
	P.~Gaucher.
	\newblock Comparing cubical and globular directed paths.
	\newblock {\em {Fund. Math.}}, 262(3):259--286, 2023.
	\newblock \href {https://doi.org/10.4064/fm219-3-2023}
	{\path{https://doi.org/10.4064/fm219-3-2023}}.
	
	\bibitem{Nonunital}
	P.~Gaucher.
	\newblock Comparing the non-unital and unital settings for directed homotopy.
	\newblock {\em Cah. Topol. G\'eom. Diff\'er. Cat\'eg.}, LXIV-2:176--197, 2023.
	
	\bibitem{ThickCubes}
	P.~Gaucher.
	\newblock Towards a theory of natural directed paths, 2023.
	\newblock \href {https://doi.org/10.48550/arXiv.2306.02792}
	{\path{https://doi.org/10.48550/arXiv.2306.02792}}.
	
	\bibitem{zbMATH07226006}
	E.~{Goubault} and S.~{Mimram}.
	\newblock {Directed homotopy in non-positively curved spaces}.
	\newblock {\em {Log. Methods Comput. Sci.}}, 16(3):55, 2020.
	\newblock Id/No 4.
	\newblock \href {https://doi.org/10.23638/LMCS-16(3:4)2020}
	{\path{https://doi.org/10.23638/LMCS-16(3:4)2020}}.
	
	\bibitem{MR1988396}
	M.~Grandis and L.~Mauri.
	\newblock Cubical sets and their site.
	\newblock {\em Theory Appl. Categ.}, 11(8):185--211, 2003.
	
	\bibitem{henry2020minimal}
	S.~Henry.
	\newblock Minimal model structures, 2020.
	\newblock \href {https://doi.org/10.48550/arXiv.2011.13408}
	{\path{https://doi.org/10.48550/arXiv.2011.13408}}.
	
	\bibitem{HKRS17}
	K.~Hess, M.~K\c{e}dziorek, E.~Riehl, and B.~Shipley.
	\newblock A necessary and sufficient condition for induced model structures.
	\newblock {\em J. Topol.}, 10(2):324--369, 2017.
	\newblock \href {https://doi.org/10.1112/topo.12011}
	{\path{https://doi.org/10.1112/topo.12011}}.
	
	\bibitem{ref_model2}
	P.~S. Hirschhorn.
	\newblock {\em Model categories and their localizations}, volume~99 of {\em
		Mathematical Surveys and Monographs}.
	\newblock American Mathematical Society, Providence, RI, 2003.
	\newblock \href {https://doi.org/10.1090/surv/099}
	{\path{https://doi.org/10.1090/surv/099}}.
	
	\bibitem{MR99h:55031}
	M.~Hovey.
	\newblock {\em Model categories}.
	\newblock American Mathematical Society, Providence, RI, 1999.
	\newblock \href {https://doi.org/10.1090/surv/063}
	{\path{https://doi.org/10.1090/surv/063}}.
	
	\bibitem{LawvereMetric}
	F.~W. Lawvere.
	\newblock Metric spaces, generalized logic, and closed categories.
	\newblock {\em Repr. Theory Appl. Categ}, 2002(1):1--37, 2002.
	
	\bibitem{coend-calculus}
	F.~Loregian.
	\newblock {\em ({Co})end calculus}, volume 468 of {\em Lond. Math. Soc. Lect.
		Note Ser.}
	\newblock Cambridge: Cambridge University Press, 2021.
	\newblock \href {https://doi.org/10.1017/9781108778657}
	{\path{https://doi.org/10.1017/9781108778657}}.
	
	\bibitem{MR1712872}
	S.~Mac~Lane.
	\newblock {\em Categories for the working mathematician}.
	\newblock Springer-Verlag, New York, second edition, 1998.
	\newblock \href {https://doi.org/10.1007/978-1-4757-4721-8}
	{\path{https://doi.org/10.1007/978-1-4757-4721-8}}.
	
	\bibitem{MR2521708}
	M.~Raussen.
	\newblock Trace spaces in a pre-cubical complex.
	\newblock {\em Topology Appl.}, 156(9):1718--1728, 2009.
	\newblock \href {https://doi.org/10.1016/j.topol.2009.02.003}
	{\path{https://doi.org/10.1016/j.topol.2009.02.003}}.
	
	\bibitem{CategoryInContext}
	E.~Riehl.
	\newblock {\em Category theory in context}.
	\newblock Aurora: Modern Math Originals. Aurora, 2016.
	
	\bibitem{MR2506258}
	J.~Rosick{\'y}.
	\newblock On combinatorial model categories.
	\newblock {\em Appl. Categ. Struct.}, 17(3):303--316, 2009.
	\newblock \href {https://doi.org/10.1007/s10485-008-9171-2}
	{\path{https://doi.org/10.1007/s10485-008-9171-2}}.
	
	\bibitem{zbMATH06722019}
	J.~Rosick{\'y}.
	\newblock Accessible model categories.
	\newblock {\em Appl. Categ. Struct.}, 25(2):187--196, 2017.
	\newblock \href {https://doi.org/10.1007/s10485-015-9419-6}
	{\path{https://doi.org/10.1007/s10485-015-9419-6}}.
	
	\bibitem{c-Reedy}
	M.~Shulman.
	\newblock Reedy categories and their generalizations, 2015.
	\newblock \href {https://doi.org/10.48550/arXiv.1507.01065}
	{\path{https://doi.org/10.48550/arXiv.1507.01065}}.
	
	\bibitem{MR4070250}
	K.~Ziemia\'{n}ski.
	\newblock Spaces of directed paths on pre-cubical sets {II}.
	\newblock {\em J. Appl. Comput. Topol.}, 4(1):45--78, 2020.
	\newblock \href {https://doi.org/10.1007/s41468-019-00040-z}
	{\path{https://doi.org/10.1007/s41468-019-00040-z}}.
	
\end{thebibliography}


\end{document}